\documentclass[12pt]{article}
\usepackage{mathrsfs}
\usepackage{amsfonts}
\usepackage{psfrag,epsfig}
\usepackage{latexsym,amsmath,amssymb}
\usepackage{exscale,relsize}
\usepackage{siunitx}
\usepackage{graphicx}
\usepackage{epstopdf}
\usepackage{indentfirst,cite}
\numberwithin{equation}{section}

\begin{document}
\title{{Ground states for fractional elliptic systems with critical growth, sublinear absorption and subcritical perturbation}
\thanks{Research supported by Natural Science Foundation of Shandong Province (No. ZR2024MA009).}
\author {Guanwei Chen$^{*}$\\
\small  \it School of Mathematical Sciences, University of Jinan, Jinan 250022, Shandong Province, China\\
\small E-mail: guanweic@163.com}}

\date{}
\maketitle

\baselineskip=8pt

\noindent \begin{abstract}
This paper studies a class of fractional elliptic systems and equations with critical exponents, sublinear absorption, subcritical perturbations, and general frequency parameters, both within a general framework involving non-local integrodifferential operators and, as a special case, for the fractional Laplacian $(-\Delta)^s$.
Our main results extend the existing literature on critical fractional elliptic problems by obtaining ground state solutions for coupled systems with sublinear absorption and a subcritical perturbation, so that the nonlinearity is not superlinear at zero; the proofs combine a Jeanjean-type monotonicity argument with a sharp mountain-pass estimate below the critical threshold, together with two concrete verifications of the energy threshold condition beyond the pure critical case.
\end{abstract}



\baselineskip=18pt

\section*{1. Introduction and Main Results}

Non-local and fractional operators of elliptic type have attracted a great
deal of attention in recent years, both for their rich theoretical structure
and for their applications to quantum mechanics, anomalous diffusion, phase
transitions, crystal dislocation, and game theory; see, e.g.,
\cite{1,5,16,18,19,33} and the references therein. A model case is the
fractional Laplacian \((-\Delta)^s\), \(s\in(0,1)\), defined (up to a
normalization constant) by
\[
-(-\Delta)^{s}u(x) = \int_{\mathbb{R}^N} \frac{u(x + y) + u(x - y) - 2u(x)}{|y|^{N+2s}} dy, \quad x \in \mathbb{R}^N.
\]

A large body of work has been devoted to fractional Schr\"odinger
equations and related non-local problems, establishing existence,
multiplicity, and qualitative properties of weak solutions. We now recall the most relevant contributions.

\textbf{$\bullet$ Subcritical fractional problems.}
In the subcritical case, the compact embeddings of the energy space lead
to abundant results. On bounded domains, many authors
\cite{20,29,31,28,B-Zhang,Zhang2017} have studied the subcritical
fractional elliptic equation
\[
\left\{
\begin{array}{ll}
-\mathcal{L}_{K_1}u - \mu u = g(x,u), & x \in \Omega, \\[4pt]
u(x) = 0, & x \in \mathbb{R}^N \setminus \Omega,
\end{array}
\right.
\]
where \(g(x,u)\) satisfies
\(\lim_{|u|\to+\infty}g(x,u)/|u|^{2^*_s-1}=0\) uniformly in \(x\in\Omega\).
As for the frequency parameter \(\mu\), \cite{20,29,B-Zhang,Zhang2017}
treat the case where \(\mu\) is sufficiently small, while \cite{28,31}
investigate the case of an arbitrary \(\mu\in\mathbb R\). The existence
of nodal and multiple solutions for the Dirichlet problem
\((-\Delta)^s u = f(x,u)\) on a bounded domain was further established
in \cite{8}. On the whole space \(\mathbb R^N\), subcritical fractional
Schr\"odinger equations have been widely studied; see, e.g.,
\cite{Felmer,Secchi,Dipierro,Byeon,Chang,36}; here the lack of
compactness comes from the translation invariance of \(\mathbb R^N\),
and not from the critical exponent. For subcritical fractional systems
and related nonlocal problems, including concave-convex systems
\cite{ChenDeng}, nonlocal H\'enon-type equations \cite{Bendahou},
perturbed fractional \(p\)-Laplacian equations \cite{38}, and
superlinear fractional \(\Phi\)-Laplacian problems \cite{Assis}, we
refer the reader to the corresponding references.

\textbf{$\bullet$ Critical fractional problems.}
For critical fractional problems, the embedding into \(L^{2^*_s}(\Omega)\)
is non-compact, making them more challenging and hence a focus of
research. The results most relevant to the present paper are the
critical fractional elliptic equations with a general kernel and the
critical fractional Laplace equations \cite{Barrios,32}. For other
related critical fractional problems, we refer the reader to fractional
\(p\)-Laplacian systems \cite{Qiao2022,Zediri2022}, fractional Kirchhoff
equations \cite{Kirchhoff2017}, and \(p\)-fractional Choquard-type
equations \cite{Sakuma2025}.

Among the works most closely related to ours are \cite{Barrios,32}. The critical fractional elliptic equation with a general kernel
\[
\left\{
\begin{array}{ll}
-\mathcal{L}_{K_1}u - \mu u = g(x,u) + |u|^{2^*_s-2}u, & x \in \Omega, \\[4pt]
u(x) = 0, & x \in \mathbb{R}^N \setminus \Omega,
\end{array}
\right.
\]
was first studied in \cite{32} for a general kernel \(K_1\) satisfying
\((K_1)\)--\((K_3)\), under the assumptions
\begin{equation}
\mu \in (0, \sigma_1),
\qquad
\lim_{|u|\to 0} \frac{g(x,u)}{u} = 0 \quad \text{uniformly in } x \in \Omega,
\tag{1.1}
\end{equation}
where \(\sigma_1\) is the first eigenvalue of \(-\mathcal{L}_{K_1}\).
The special case \(\mathcal{L}_{K_1}=(-\Delta)^s\) was treated in
\cite{Barrios}. The limitations of these two works are summarized
below.

For the existing critical fractional problems \cite{32,Barrios}, the
frequency parameter is required to satisfy \(0<\mu<\sigma_1\), the
nonlinearity is required to be superlinear at zero in the sense of
(1.1), and the concrete verification of the energy threshold condition
is carried out only in the pure critical case. None of these works, however, covers the combination of a sublinear
absorption term with a subcritical perturbation, a general frequency
range, a threshold verification in their presence, and ground state
solutions.

Motivated by the above results, we study the following critical fractional
elliptic system with
general kernels and a sublinear absorption term:
$$
\left\{\begin{array}{ll}
-\mathcal{L}_{K} U - M U = -\gamma(x)|U|^{\delta-2}U+\nabla F(x, U)+|U|^{2^*_s-2}U , & x \in \Omega, \\[4pt]
U = \mathbf{0}, & x \in \mathbb{R}^N \setminus \Omega,\\[4pt]
\end{array}\right.
\eqno(1.2)
$$
where \(U=(u_1(x),\cdots,u_k(x))\) is a \(k\)-dimensional vector function,
\(\Omega\subset\mathbb R^N\) is a bounded smooth domain, \(N>2s\),
\(s\in(0,1)\), \(2^*_s:=\frac{2N}{N-2s}\), and
\(M=\operatorname{diag}(\mu_1,\dots,\mu_k)\) is the frequency matrix. The sublinear absorption term is defined by
$$
|U|^{\delta-2}U:=
\left\{\begin{array}{ll}
0, &|U|=0,  \\[4pt]
\frac{U}{|U|^{2-\delta}},  &|U|\neq0,
\end{array}\right.
\quad \delta\in(1,2).
$$
Here \(\nabla F(x, U) = (f_1,\dots,f_k)\) is a subcritical perturbation, with
\(F(x, U) \in C^1(\Omega \times \mathbb{R}^k, \mathbb{R})\), and the
operator \(\mathcal{L}_{K} U := (\mathcal{L}_{K_1} u_1, \dots,
\mathcal{L}_{K_k} u_k)\) consists of the integro-differential operators
$$
\mathcal{L}_{K_i} u_i(x):= \int_{\mathbb{R}^N} \left[ u_i(x + y) + u_i(x - y) - 2u_i(x) \right] K_i(y) dy.
$$
The kernels \(K_i\) (\(i=1,\dots,k\)) satisfy:

\smallskip
\noindent${\bf (K_1)} \ \ K_i : \mathbb{R}^N \setminus \{0\} \rightarrow (0, +\infty)$,  and $\min\{|x|^2, 1\} K_i \in L^1(\mathbb{R}^N,\mathbb{R})$;\\
${\bf (K_2)}$ There exist $\rho> 0$ and $s \in (0, 1)$ such that $K_i(x) \geq \rho |x|^{-(N+2s)}$ for all $x \neq 0$;\\
${\bf (K_3)} \ \ K_i(x) = K_i(-x)$ for all $x \neq 0.$

\smallskip
Under \((K_1)\)--\((K_3)\), Proposition 9 in \cite{31} implies that the
spectrum of \(-\mathcal{L}_{K_i}\) consists entirely of positive
eigenvalues
$$
0<\sigma^i_1<\sigma^i_2\le\sigma^i_3\le\sigma^i_m\le\cdots\to\infty,\quad i=1,\cdots,k.
\eqno(1.3)
$$
For refined eigenvalue bounds on the Dirichlet fractional Laplacian we
refer to \cite{39}, and for the second eigenvalue of the fractional
\(p\)-Laplacian to \cite{Brasco}.

For the particular choice \(K_i(x)\equiv |x|^{-(N+2s)}\), the system (1.2)
reduces, up to a normalization constant, to the fractional Schr\"odinger
system
$$
\left\{\begin{array}{ll}
(-\Delta)^{s} U - M U = -\gamma(x)|U|^{\delta-2}U+\nabla F(x, U)+|U|^{2^*_s-2}U, & x \in \Omega, \\[4pt]
U = \mathbf{0}, & x \in \mathbb{R}^N \setminus \Omega,
\end{array}\right.
\eqno(1.4)
$$
where \((-\Delta)^s U:=((-\Delta)^{s} u_1,\cdots,(-\Delta)^{s}u_k)\).
System (1.4) corresponds to standing waves
\(\psi(x,t)=e^{\ell Mt}U(x)\) of a time-dependent fractional
Schr\"odinger system; \(U\) is the amplitude and \(M\) the frequency
matrix.

\medskip
\noindent\textbf{Definition 1.1 (Function spaces).}
For \(i = 1,\dots,k\), let \(E_i\) be the linear space of Lebesgue
measurable functions \(u : \mathbb{R}^N \to \mathbb{R}\) such that
\(u|_\Omega \in L^2(\Omega)\) and
\[
(x,y) \mapsto (u(x)-u(y))\sqrt{K_i(x-y)} \in L^2\!\Bigl(\mathbb{R}^{2N} \setminus \bigl[(\mathbb{R}^N \setminus \Omega) \times (\mathbb{R}^N \setminus \Omega)\bigr], dx\,dy\Bigr).
\]
Define
\[
\widetilde{E}_i := \{ u \in E_i : u = 0 \text{ a.e. in } \mathbb{R}^N \setminus \Omega \}, \qquad
X := \widetilde{E}_1 \times \cdots \times \widetilde{E}_k,
\]
which are non-empty since \(C_0^2(\Omega) \subset \widetilde{E}_i \subset E_i\) by Lemma 5.1 in \cite{30}.

\medskip
\noindent\textbf{Definition 1.2 (Weak solutions).}
From \cite{28,29,31}, \((\widetilde{E}_i, \|\cdot\|_{\widetilde{E}_i})\) are Hilbert spaces
 with norms $\|u_i\|^2_{\widetilde{E}_i}:=(u_i,u_i)_{\widetilde{E}_i}$ and inner products
$$
(u_i,\varphi_i)_{\widetilde{E}_i}:=\int_{Q} (u_i(x) - u_i(y))(\varphi_i(x) - \varphi_i(y)) K_i(x-y) dxdy,
\eqno(1.5)
$$
where \(Q := \mathbb{R}^{2N} \setminus [(\mathbb{R}^N \setminus \Omega) \times (\mathbb{R}^N \setminus \Omega)]\).
A function \(U = (u_1,\dots,u_k) \in X\) is a weak solution of (1.2) if
for every \(\varphi = (\varphi_1,\dots,\varphi_k) \in X\),
$$
\begin{array}{rcl}
\mathlarger{\sum_{i=1}^{k}}(u_i,\varphi_i)_{\widetilde{E}_i}
&=& \mathlarger{\sum_{i=1}^{k}}\mu_i(u_i,\varphi_i)_{L^2}-\int_{\Omega} \gamma(x)|U|^{\delta-2}U \cdot\varphi dx
+ \int_\Omega |U|^{2^*_s-2}U\cdot\varphi  dx\\
\\
&&+ \int_{\Omega} \nabla F(x, U)\cdot\varphi dx.
\end{array}
\eqno(1.6)
$$

\medskip
We impose the following assumptions on the weight function \(\gamma\),
the nonlinearity \(F\), and the frequency matrix \(M\): 

\smallskip
\noindent${\bf (\Gamma_{1})} \ \ \gamma(x)\in L^\infty(\Omega)$ and $\gamma(x)\ge0$ a.e. in $\Omega$;\\
${\bf (F_{1})} \ \  F\in C^1,\ F(x,U)\ge0,\ F(x,\mathbf{0})=0$, and $\underset{x \in \Omega, \|U\| \le C}{\text{ess sup}}|\nabla F(x,U)| < +\infty$, $\forall C > 0$;\\
${\bf (F_{2})}  \ \  \lim\limits_{|U|\to0}\frac{|\nabla F(x,U)|}{|U|}=0$ uniformly in $x\in\Omega$;\\
${\bf (F_{3})} \ \  \lim\limits_{|U|\to+\infty}\frac{|\nabla F(x,U)|}{|U|^{2^*_s-1}}=0$ uniformly in $x \in \Omega$;\\
${\bf (F_{4})} \ \   \widetilde F(x,U):=\tfrac12\nabla F(x,U)\cdot U-F(x,U)\ge0$,
a.e. $x\in\Omega$, $\forall U\in\mathbb R^k\setminus\{\mathbf{0}\}$;\\ 
${\bf (M_1)} \ \   -\infty<\mu_i<\sigma_1^i,\ i=1,\cdots,k.$

\smallskip

Under the above assumptions, the functional 
$$\begin{array}{rcl}
\Phi(U) &=& \frac{1}{2} \|U\|^2 -\frac{1}{2}\sum_{i=1}^{k}\mu_i\|u_i\|_{L^2}^2+\frac{1}{\delta}\int_\Omega \gamma(x)|U|^{\delta} dx\\
\\
&&-\frac{1}{2^*_s}\|U\|_{2^*_s}^{2^*_s}-\int_\Omega F(x, U) dx , \quad \forall U \in X,
\end{array}
\eqno(1.7)
$$
is well defined and of class $C^1$ on $X$, and weak solutions of the system (1.2) are precisely its critical points.
Here \(\|\cdot\|_{2^*_s}\) denotes the norm in \([L^{2^*_s}(\Omega)]^k\) and
\begin{equation}
\|U\|^2 := \sum_{i=1}^{k}\|u_i\|^2_{\widetilde{E}_i}. \tag{1.8}
\end{equation}
The best Sobolev constant is
\begin{equation}
S_K := \inf_{V \in X \setminus \{\mathbf{0}\}}
\frac{\|V\|^2}{\|V\|_{2^*_s}^2}.
\tag{1.9}
\end{equation}
This is the vector-valued analogue of the constant \(S_{K_1}\) in
\cite[formula (1.17)]{32}; equivalently, it can be characterized by the
scalar constants \(S_{K_i}\) associated with each component. Since the
embedding \(X\hookrightarrow[L^{2^*_s}(\Omega)]^k\) is continuous
(Lemma 2.1), we have \(S_K>0\).

Due to the critical exponent, the Palais--Smale (PS)
condition for \(\Phi\) holds only below a threshold determined by
\(S_K\); we thus impose the following energy threshold condition.

\smallskip
\noindent${\bf (S_{1})} \quad$ There exists $U_0\in X\setminus\{\mathbf{0}\}$ such that
$$
\sup_{\zeta\ge0}\Phi(\zeta U_0)
<
2^{-\frac{N-2s}{2s}}\frac{s}{N}S_K^{N/(2s)}.
$$

\smallskip
\noindent\textbf{Theorem 1.1 (Abstract result).} {\it
Assume \((K_1)\)--\((K_3)\), \((\Gamma_1)\), \((F_1)\)--\((F_4)\), \((M_1)\)
and \((S_1)\). Then $(1.2)$ has a ground state solution
\(\widetilde U\in X\setminus\{0\}\) such that}
\[
\Phi(\widetilde U)=\widetilde c
:=\inf\{\Phi(U):U\in X\setminus\{0\},\ \Phi'(U)=0\}>0.
\]

\smallskip
Next, we verify the energy threshold condition \((S_1)\) under different assumptions, and Theorem 1.1
then gives a ground state solution of the systems (1.2) and (1.4).

\smallskip
\noindent\textbf{Notation 1.1.}
Throughout, \(u_0\) denotes the positive first eigenfunction of
\(-\mathcal L_{K_1}\) associated with \(\sigma_1^1\), normalized by
\(\|u_0\|_{L^2}=1\), and set
\(U_0:=(u_0,0,\dots,0)\in X\).

\smallskip
\noindent\textbf{Theorem 1.2 (Small-perturbation verification of \((S_1)\) for (1.2)).}
{\it Assume \((K_1)\)--\((K_3)\) and \((M_1)\).

\noindent\textbf{(i) The pure critical case.}
Suppose \(\gamma\equiv0\), \(F\equiv0\) in \(\Omega\times\mathbb R^k\),
and \(\mu_1\) as in \((M_1)\) satisfies
\[
\mu_1>\sigma_1^1-2^{-\frac{N-2s}{N}}S_K\|u_0\|_{L^{2^*_s}}^2.
\tag{1.10}
\]
Then \((S_1)\) holds with \(U_0=(u_0,0,\dots,0)\in X\).

\noindent\textbf{(ii) The small-perturbation case for (1.2).}
Suppose \((\Gamma_1)\), \((F_1)\)--\((F_4)\) and \((1.10)\) hold.
Moreover, assume that there exists a sufficiently small
\(\varepsilon_0=\varepsilon_0(U_0)>0\) such that
\[
\|\gamma\|_{L^\infty(\Omega)}\le\varepsilon_0,
\qquad
|\nabla F(x,U)|\le\varepsilon_0|U|+\varepsilon_0|U|^{2^*_s-1},
\quad \forall(x,U)\in\Omega\times\mathbb R^k.
\tag{1.11}
\]
Then \((S_1)\) holds with \(U_0=(u_0,0,\dots,0)\in X\).

Consequently, in both cases, Theorem $1.1$ yields a ground state solution
\(\widetilde U\in X\setminus\{0\}\) of \((1.2)\) with
\(\Phi(\widetilde U)=\widetilde c>0\).}

\smallskip
\noindent\textbf{Theorem 1.3 (Quantitative absorption verification of \((S_1)\) for (1.2)).}
{\it Assume \((K_1)\)--\((K_3)\), \((\Gamma_1)\), \((F_1)\)--\((F_4)\), \((M_1)\) and $(1.10)$, and let \(u_0\), \(U_0\) be as in the Notation above. Then there exists \(\varepsilon_1=\varepsilon_1(N,s,\sigma_1^1,
\mu_1,\delta,S_K,u_0)>0\) such that if
\[
\int_\Omega\gamma(x)|u_0|^\delta\,dx<\varepsilon_1,
\]
then \((S_1)\) holds with \(U_0\). Consequently, by Theorem $1.1$, the
system $(1.2)$ admits a ground state solution
\(\widetilde U\in X\setminus\{0\}\) with
\(\Phi(\widetilde U)=\widetilde c>0\).}

\smallskip
\noindent\textbf{Remark 1.1.}
Theorem 1.2 (i) treats the pure critical case \(\gamma\equiv0\) and
\(F\equiv0\). Theorem 1.2 (ii) requires \(\gamma\) and \(\nabla F\) to be
small in \(L^\infty\). Theorem 1.3 weakens these assumptions: it only
requires the weighted integral
\(\int_\Omega\gamma|u_0|^\delta\) to be small, and imposes no
smallness on \(F\). Thus Theorem 1.3 contains Theorem 1.2 (ii) as a
special case and extends the verification of \((S_1)\) to a strictly
larger class of nonlinearities. Theorem 1.3 provides a quantitative sufficient condition for
\((S_1)\): the weighted integral \(\int_\Omega\gamma|u_0|^\delta\)
must be small compared with the gap between the best Sobolev energy
threshold and the maximal energy along the ray generated by \(U_0\);
no smallness of \(F\) in \(L^\infty\) is required, and only this single
weighted quantity of the absorption coefficient enters the condition.

\smallskip
\noindent\textbf{Theorem 1.4 (Results of the system (1.4)).}
{\it Assume \((M_1)\) holds with \(\sigma_1^i\) replaced by
\(\sigma_{1,s}\), \(i=1,\dots,k\), where \(\sigma_{1,s}\) is the first
eigenvalue of \((-\Delta)^s\) with homogeneous Dirichlet boundary
condition. Let \(u_{0,s}\) be the corresponding positive first
eigenfunction normalized by \(\|u_{0,s}\|_{L^2}=1\), and let \(S_s\) be
the best fractional Sobolev constant in
\(H^s(\mathbb R^N)\hookrightarrow L^{2^*_s}(\mathbb R^N)\). Set
\(U_0:=(u_{0,s},0,\dots,0)\in X\).

\noindent\textbf{(i) The pure critical case.}
Suppose \(\gamma\equiv0\), \(F\equiv0\), and \(\mu_1\) satisfies
\[
\mu_1>\sigma_{1,s}-2^{-\frac{N-2s}{N}}S_s\|u_{0,s}\|_{2^*_s}^2.
\tag{1.12}
\]
Then \((S_1)\) holds with \(U_0\).

\noindent\textbf{(ii) The small-perturbation case.}
Suppose \((\Gamma_1)\), \((F_1)\)--\((F_4)\), \((1.11)\) and \((1.12)\)
hold. Then \((S_1)\) holds with \(U_0\).

\noindent\textbf{(iii) The quantitative absorption case.}
Suppose \((\Gamma_1)\), \((F_1)\)--\((F_4)\) and \((1.12)\) hold. Then
there exists
\(\varepsilon_{1,s}
=\varepsilon_{1,s}(N,s,\sigma_{1,s},\mu_1,\delta,S_s,u_{0,s})>0\)
such that if
\[
\int_\Omega\gamma(x)|u_{0,s}|^\delta\,dx<\varepsilon_{1,s},
\]
then \((S_1)\) holds with \(U_0\).

Consequently, in all three cases, Theorem $1.1$ applied to the special
case \(K_i(x)=|x|^{-(N+2s)}\) implies that the system \((1.4)\) admits
a ground state solution \(\widetilde U\in X\setminus\{0\}\) satisfying
\(\Phi(\widetilde U)=\widetilde c>0\).}

\smallskip
\noindent\textbf{Remark 1.2.}
The main difficulties and our strategy for overcoming them are as follows.

\noindent\textbf{(1)} {\it Lack of compactness.} Due to the critical
exponent \(2^*_s\), the embedding
\(X\hookrightarrow[L^{2^*_s}(\Omega)]^k\) is not compact, so the PS
condition for \(\Phi\) does not hold globally. We impose and
verify the energy threshold condition \((S_1)\), and estimate the
Mountain-Pass level precisely.

\noindent\textbf{(2)} {\it Failure of superlinearity at zero.} The
sublinear absorption term \(-\gamma(x)|U|^{\delta-2}U\) with
\(\delta\in(1,2)\) destroys the superlinearity at zero required by the
standard Mountain-Pass argument used in \cite{32,Barrios}. We use the method of Jeanjean \cite{Jeanjean}, introducing approximating
functionals \(\Phi_\lambda=A-\lambda B\), \(\lambda\in[1,2]\), and
passing to the limit \(\lambda\to1\).

\noindent\textbf{(3)} {\it Difficulty in verifying \((S_1)\).} Because of
the sublinear absorption and the subcritical perturbation
\(\nabla F(x,U)\), the condition \((S_1)\) cannot be verified by the
same computations as in the pure critical case \cite{32,Barrios}. We provide
two types of concrete verification: a small-perturbation one in
Theorem 1.2, and a quantitative one in Theorem 1.3 requiring only a
single weighted quantity of the absorption coefficient to be small and
imposing no smallness on \(F\). These verifications are then realized
for the concrete system (1.4) in Theorem 1.4.

\medskip

\noindent\textbf{Remark 1.3.}
We extend the existing critical fractional results
most closely related to our setting \cite{32,Barrios} in five
directions:
\\
\textbf{(1)} a general frequency range \(\mu_i<\sigma_1^i\),
which includes negative values and is strictly wider than the smallness
condition used in \cite{32,Barrios};
\\
\textbf{(2)} a sublinear absorption term \(-\gamma(x)|U|^{\delta-2}U\)
with \(\delta\in(1,2)\), so that the overall nonlinearity is no longer
superlinear at zero in the sense of (1.1), in contrast with
\cite{32,Barrios};
\\
\textbf{(3)} a concrete and quantitative verification of \((S_1)\) in the
presence of both the sublinear absorption and the subcritical
perturbation (Theorems 1.2--1.3), together with its concrete
realization for system (1.4) (Theorem 1.4); in particular, the
quantitative condition in Theorem 1.3 requires only the weighted
integral \(\int_\Omega\gamma|u_0|^\delta\) to be small and imposes no
smallness on \(F\), whereas the existing verifications in
\cite{32,Barrios} are restricted to the pure critical case;
\\
\textbf{(4)} coupled systems rather than scalar equations in \cite{32,Barrios};
\\
\textbf{(5)} ground state solutions rather than merely nontrivial ones in \cite{32,Barrios}.

\smallskip
\noindent\textbf{Example 1.1.}
A representative nonlinearity satisfying \((F_1)\)--\((F_4)\) is
\[
F(x,U):=a(x)|U|^q+b(x)|U|^r,\qquad 2<q<r<2^*_s,\qquad
0\le a,b\in L^\infty(\Omega).
\]
Indeed, \((F_1)\) follows from \(q,r>1\) and \(a,b\ge0\); \((F_2)\) and
\((F_3)\) follow from \(2<q<r<2^*_s\). Moreover,
\[
\widetilde F(x,U)
=\frac{q-2}{2}a(x)|U|^q+\frac{r-2}{2}b(x)|U|^r\ge0,
\]
so \((F_4)\) also holds.
If, in addition, \(\|a\|_{L^\infty}\) and \(\|b\|_{L^\infty}\) are
sufficiently small, then (1.11) also holds. Indeed, splitting the
estimate according to \(|U|\le1\) and \(|U|>1\), we have
\[
|\nabla F(x,U)|
\le\bigl(q\|a\|_{L^\infty}+r\|b\|_{L^\infty}\bigr)
\bigl(|U|+|U|^{2^*_s-1}\bigr),
\]
so (1.11) holds with
\(\varepsilon_0=q\|a\|_{L^\infty}+r\|b\|_{L^\infty},\)
which can be made arbitrarily small by choosing
\(\|a\|_{L^\infty},\|b\|_{L^\infty}\) sufficiently small.

Thus, for all \(a,b\ge0\), this \(F\) satisfies \((F_1)\)--\((F_4)\) and
hence is admissible in Theorem 1.3; if in addition
\(\|a\|_{L^\infty},\|b\|_{L^\infty}\) are sufficiently small, then
(1.11) holds and the example is also admissible in Theorem 1.2(ii).

\smallskip
\smallskip
The paper is organized as follows. Section 2 introduces the variational
framework and preliminary lemmas. Section 3 contains the proofs of our
main results.

\section*{2. Variational Framework and Preliminary Lemmas}
In the following, we will use \(``\to"\) to denote strong convergence and \(``\rightharpoonup"\) to denote weak convergence.

From \cite{28,29,31}, we know \((X, \|\cdot\|)\) is a Hilbert space with the inner product
$$
(U, \varphi) =\mathlarger{\sum_{i=1}^{k}}(u_i,\varphi_i)_{\widetilde{E}_i},\quad U=(u_1,\cdots,u_k), \ \varphi=(\varphi_1,\cdots,\varphi_k) \in X.
\eqno(2.1)
$$
Under our assumptions, the functional \(\Phi\in C^1(X, \mathbb{R})\) and
$$\begin{array}{rcl}
\langle \Phi'(U), \varphi \rangle &=& (U, \varphi)-\mathlarger{\sum_{i=1}^{k}}\mu_i(u_i,\varphi_i)_{L^2}
+ \int_\Omega \gamma(x) |U|^{\delta-2}U\cdot\varphi  dx\\
\\
&& - \int_\Omega |U|^{2^*_s-2}U\cdot\varphi  dx-  \int_\Omega \nabla F(x, U(x))\cdot\varphi dx,
\end{array}\eqno(2.2)
$$
for all \(U=(u_1,\cdots,u_k)\) and \(\varphi=(\varphi_1,\cdots,\varphi_k) \in X\). Thus, the problem of finding a nontrivial (weak) solution of the system (1.2) has been reduced to that of seeking a nonzero critical point of the functional \(\Phi\) on \(X\).

\medskip
\noindent\textbf{Lemma 2.1 (Embedding theorem).} If \((K_1)\)--\((K_3)\) hold, then the following embeddings are continuous:
$$
(X,\|\cdot\|) \hookrightarrow [L^t(\Omega)]^k:=L^t\times\cdots\times L^t,\quad \forall t \in [1, 2^*_s],
$$
i.e., there exists \(\tau_t > 0\) such that
$\|U\|_t \leq \tau_t \|U\|$ for all $U \in X$ and all $t \in [1, 2^*_s],$
where \(\|\cdot\|_t\) denote the norms in \([L^t(\Omega)]^k\). Moreover, the following embeddings are compact:
$$
(X,\|\cdot\|) \hookrightarrow [L^l(\Omega)]^k,\quad \forall l \in [1, 2^*_s),
$$
\textbf{Proof.} The scalar case is proved in \cite[Lemma 8]{29}: for each
\(i=1,\dots,k\), the embedding
\(\widetilde E_i\hookrightarrow L^t(\Omega)\) is continuous for every
\(t\in[1,2^*_s]\), and compact for every \(t\in[1,2^*_s)\).

For the vector-valued case, let \(U=(u_1,\dots,u_k)\in X\). Then
\[
\|U\|_t^2
=\sum_{i=1}^k\|u_i\|_{L^t(\Omega)}^2
\le \sum_{i=1}^k\bigl(\tau_t^{(i)}\bigr)^2\|u_i\|_{\widetilde E_i}^2
\le \Bigl(\max_{i=1,\dots,k}\tau_t^{(i)}\Bigr)^2\|U\|^2,
\]
where \(\tau_t^{(i)}\) is the embedding constant for the \(i\)-th
component. Hence the embedding \(X\hookrightarrow[L^t(\Omega)]^k\) is
continuous for every \(t\in[1,2^*_s]\).

If \(\{U_n\}\) is bounded in \(X\), then each component sequence
\(\{u_n^i\}\) is bounded in \(\widetilde E_i\). By the compactness of
the scalar embedding \(\widetilde E_i\hookrightarrow L^l(\Omega)\),
\(l\in[1,2^*_s)\), each component has a strongly convergent subsequence
in \(L^l(\Omega)\). A standard diagonal argument then yields a
subsequence of \(\{U_n\}\) converging strongly in \([L^l(\Omega)]^k\).
Hence \(X\hookrightarrow[L^l(\Omega)]^k\) is compact for every
\(l\in[1,2^*_s)\).$\quad\Box$

\medskip
\noindent\textbf{Lemma 2.2 (Inequalities involving linear terms and norms).}
{\it If \((K_1)\)--\((K_3)\) and \((M_1)\) hold, then we have
$$
\|U\|^2-\mathlarger{\sum_{i=1}^{k}}\mu_i\|u_i\|_{L^2}^2\ge\overline{D}\|U\|^2,\quad \forall U\in X,
\eqno(2.3)$$
where
$\overline{D}:=
\min\left\{1, \min_{i=1,\cdots,k}\left(1-\frac{\mu_i}{\sigma_1^i}\right)\right\}>0.$}\\
\textbf{Proof.} If \(0<\mu_i<\sigma_1^i\), \(i=1,\cdots,k\), by (1.8) and (1.3), we have
$$\begin{array}{rcl}
\|U\|^2-\mathlarger{\sum_{i=1}^{k}}
\left(1-\frac{\mu_i}{\sigma_1^i}\right)(u_i,u_i)_{\widetilde{E}_i}
&=\mathlarger{\sum_{i=1}^{k}}(u_i,u_i)_{\widetilde{E}_i}-\mathlarger{\sum_{i=1}^{k}}\left(1-\frac{\mu_i}{\sigma_1^i}\right)(u_i,u_i)_{\widetilde{E}_i}\\
\\
&=\mathlarger{\sum_{i=1}^{k}}\frac{\mu_i}{\sigma_1^i}(u_i,u_i)_{\widetilde{E}_i}
\ge\mathlarger{\sum_{i=1}^{k}}\mu_i\|u_i\|_{L^2}^2, \quad \forall U\in X.
\end{array}
$$
So we have
$$\begin{array}{rcl}
\|U\|^2-\mathlarger{\sum_{i=1}^{k}}\mu_i\|u_i\|_{L^2}^2
&\ge&\mathlarger{\sum_{i=1}^{k}}\left(1-\frac{\mu_i}{\sigma_1^i}\right)(u_i,u_i)_{\widetilde{E}_i} \\
\\
&\ge&\mathlarger{\min_{i=1,\cdots,k}}\left(1-\frac{\mu_i}{\sigma_1^i}\right)\|U\|^2, \quad \forall U\in X.
\end{array}$$
If \(-\infty<\mu_i\le0\), \(i=1,\cdots,k\), we have
$$
\|U\|^2-\mathlarger{\sum_{i=1}^{k}}\mu_i\|u_i\|_{L^2}^2
\ge\|U\|^2,\quad \forall U\in X.
$$
Therefore, Lemma 2.2 holds.$\quad\Box$

\medskip
Our method, inspired by Jeanjean (see \cite{Jeanjean}), departs from the standard PS sequence for the true functional \(\Phi\) defined in (1.7). Instead, we initiate our analysis with a sequence of exact critical points derived from the following neighboring functionals \(\Phi_\lambda\in C^1(X, \mathbb{R})\), where \(\Phi_\lambda\) defined as follows:
$$\begin{array}{rcl}
\Phi_\lambda(U)&:=&A(U)-\lambda B(U)\\
\\
&:=&\frac{1}{2}\|U\|^2
- \frac{1}{2}\mathlarger{\sum_{i=1}^{k}}\mu_i\|u_i\|_{L^2}^2+\frac{1}{\delta}\int_\Omega \gamma(x)|U|^{\delta} dx\\
\\
&&-\lambda\left(\frac{1}{2^*_s}\|U\|_{2^*_s}^{2^*_s}+\int_\Omega F(x, U) dx\right)  , \quad \forall \lambda\in J:=[1,2], \quad \forall U\in X.
\end{array}
\eqno(2.4)
$$

\medskip
\noindent\textbf{Lemma 2.3 (Mountain-pass geometry).} {\it Under the conditions of Theorem $1.1$, we have \\
$(i)$ $\Phi_\lambda(\overline{U})\leq0$ for some $\overline{U}\in X\backslash\{\mathbf{0}\}$, $\forall \lambda\in[1,2]$. \\
$(ii)$ $\Phi_{\lambda}(\widetilde{U})>0$ for some $\widetilde{U}\in X\backslash\{\mathbf{0}\}$ with $\|\widetilde{U}\|<\|\overline{U}\|$, $\forall \lambda\in[1,2]$.}\\
\textbf{Proof.} By $(F_{1})$--$(F_{3})$ and
$$
F(x,U)=\int_0^{1}\frac{\nabla F(x,tU)\cdot tU}{t}dt,
$$
we have that for any \(\varepsilon>0\) there exists \(C_\varepsilon>0\) such that
$$
|\nabla F(x,U)|\leq 2^{*}_s \varepsilon |U|^{2^{*}_s-1} + C_\varepsilon, \quad \forall(x,U)\in\Omega\times\mathbb{R}^k,
\eqno(2.5)
$$
$$
|F(x,U)|\leq \varepsilon |U|^{2^{*}_s}+ C_{\varepsilon} |U|, \quad \forall(x,U)\in\Omega\times\mathbb{R}^k.
\eqno(2.6)
$$
Similarly, for any \(\varepsilon>0\) there exists \(\delta_\varepsilon>0\) such that
$$
|\nabla F(x,U)|\leq 2\varepsilon|U| + 2^{*}_s\delta_\varepsilon|U|^{2^{*}_s-1}, \quad \forall(x,U)\in\Omega\times\mathbb{R}^k,
\eqno(2.7)
$$
$$
|F(x, U)| \leq \varepsilon |U|^{2} + \delta_\varepsilon|U|^{2^{*}_s}, \quad \forall(x,U)\in\Omega\times\mathbb{R}^k.
\eqno(2.8)
$$

{\bf Part 1.} We prove \((i)\) holds.
For the \(U_0\) given in \((S_1)\) and \(l>0\), by (2.6) (taking \(\varepsilon<\frac{1}{2^*_s}\)), Lemma 2.1, \(\lambda\in [1,2]\) and \(2^*_s>2>\delta>1\), we have that
$$\begin{array}{cl}
\Phi_{\lambda}(lU_0)&=\frac{1}{2}l^2\|U_0\|^2-
 \frac{1}{2}\mathlarger{\sum_{i=1}^{k}}\mu_il^2\|u_{0,i}\|_{L^2}^2+\frac{1}{\delta}l^{\delta}\int_\Omega \gamma(x)|U_0|^{\delta} dx\\
 \\
&\quad-\lambda\frac{l^{2^*_s}}{2^*_s}\|U_0\|_{2^*_s}^{2^*_s}-\lambda\int_\Omega F(x, lU_0) dx\\
\\
&\le-l^{2^*_s}\left(\frac{1}{2^*_s}-\varepsilon\right)\|U_0\|_{2^*_s}^{2^*_s} +l^2\|U_0\|^2 +l^2\mathlarger{\max_{i=1,\cdots,k}}|\mu_i|\|U_0\|_2^2\\
\\
&\quad +l^\delta\|\gamma\|_{L^\infty}\|U_0\|_{\delta}^{\delta} +lC_{\varepsilon}\|U_0\|_1\\
\\
&\to-\infty,\quad l\to+\infty.
\end{array}
\eqno(2.9)
$$
Therefore, by choosing \(\overline{U}:=l_0U_0\) with \(l_0>0\) sufficiently large, we have
\[
\Phi_{\lambda}(\overline{U})<0.
\]

{\bf Part 2.} We prove \((ii)\) holds.
By the fact that \(\gamma(x)\ge0\) for all \(x\in\Omega\), Lemma 2.2 and (2.8), we have there is a small enough \(r>0\) such that
$$\begin{array}{cl}
\Phi_{\lambda}(rV)&=\frac{1}{2}r^2\|V\|^2-
 \frac{1}{2}\mathlarger{\sum_{i=1}^{k}}\mu_ir^2\|v_i\|_{L^2}^2+\frac{1}{\delta}r^{\delta}\int_\Omega \gamma(x)|V|^{\delta} dx\\
 \\
&\quad-\lambda\frac{r^{2^*_s}}{2^*_s}\|V\|_{2^*_s}^{2^*_s}-\lambda\int_\Omega F(x, rV) dx\\
\\
&\ge \frac{1}{2}\overline{D}\|rV\|^2
-r^{2^*_s}\frac{\lambda}{2^*_s}\tau_{2^*_s}^{2^*_s}\|V\|^{2^*_s}-r^2\varepsilon \lambda\tau^{2}_2\|V\|^{2}
- r^{2^*_s}\lambda\delta_{\varepsilon} \tau^{2^*_s}_{2^*_s}\|V\|^{2^*_s}\\
\\
&= r^2\left(\frac{1}{2}\overline{D}\|V\|^2-\varepsilon \lambda\tau^{2}_2\|V\|^{2}\right)
-r^{2^*_s}\lambda\tau_{2^*_s}^{2^*_s}\|V\|^{2^*_s}\left(\frac{1}{2^*_s}
-\delta_{\varepsilon} \right)\\
\\
&>0\quad (\varepsilon \ \mbox{is sufficiently small}),\quad \forall \ V\in X\backslash\{\mathbf{0}\}.
\end{array}\eqno(2.10)
$$
Therefore, \((ii)\) holds by choosing \(\widetilde{U}:=\widetilde{r_0}V\) with \(\widetilde{r_0}>0\) sufficiently small and \(\|\widetilde{U}\|<\|\overline{U}\|\).$\quad\Box$

\medskip
\noindent\textbf{Lemma 2.4 (Bounded PS sequence).} {\it Under the conditions of Theorem $1.1$, for almost every $\lambda\in[1,2]$, there is a bounded sequence $\{U_j\}\subset X$ satisfying}
$$
\lim_{j\to\infty}\Phi_\lambda(U_j)= c_\lambda>0, \quad
\lim_{j\to\infty}\Phi'_\lambda(U_j)=0.
\eqno(2.11)$$
\textbf{Proof.} Let
$$
\Gamma:=\{\gamma\in C([0,1], X):\ \gamma(0)=\mathbf{0}\ \mbox{and}\
\gamma(1)=\overline{U}\},
\eqno(2.12)
$$
where \(\overline{U}=l_0U_0\) is given in the proof of Part 1 in Lemma 2.3, then \(\Phi_\lambda(\overline{U})<0\) in Lemma 2.3 and \(\Phi_{\lambda}(\mathbf{0})=0\) by \(F(x,\mathbf{0})\equiv0\) imply that
$$
c_\lambda:=\inf_{\gamma\in\Gamma}\sup_{t\in[0,1]}\Phi_{\lambda}(\gamma(t))>0=
\max\{\Phi_{\lambda}(0),\Phi_{\lambda}(\overline{U})\},\quad \forall\lambda\in J=[1,2].
\eqno(2.13)$$
It implies that \(\Phi_\lambda\) has a Mountain-Pass geometry.

On the other hand, by H\"{o}lder's inequality, Lemmas 2.1 and 2.2, and \(2>\delta>1\), we have
$$\begin{array}{rcl}
A(U)&=&\frac{1}{2}\|U\|^2
- \frac{1}{2}\mathlarger{\sum_{i=1}^{k}}\mu_i\|u_i\|_{L^2}^2
+\frac{1}{\delta}\int_\Omega \gamma(x)|U|^{\delta} dx  \\
\\
&\ge&\frac{1}{2}\overline{D}\|U\|^2-\|\gamma\|_{L^\infty}\tau_\delta^\delta\|U\|^\delta\\
\\
&\rightarrow&+\infty,\quad
\|U\|\rightarrow+\infty,
\end{array}
$$
where \(\overline{D}>0\) is given in Lemma 2.2. Moreover, \((F_1)\) implies
$$
B(U)=\frac{1}{2^*_s}\|U\|_{2^*_s}^{2^*_s}+\int_\Omega F(x, U) dx\ge0.
$$

Therefore, the proof is finished by Theorem 1.1 in \cite{Jeanjean} (see Appendix).$\quad\Box$

\medskip
\noindent\textbf{Lemma 2.5 (Existence of a nontrivial critical point).} {\it Under the conditions of Theorem $1.1$, for any fixed
\(\lambda\in[1,2]\), if \(\{U_j\}\subset X\) is bounded and
satisfies
$$
0<\lim_{j\to\infty}\Phi_\lambda(U_j)= c_\lambda, \quad
\lim_{j\to\infty}\Phi'_\lambda(U_j)=0,
\eqno(2.14)
$$
then after passing to a subsequence (still denoted by
$\{U_j\}$),}
$$
U_j\rightharpoonup U_\lambda\neq0 \ in \ X, \quad
\Phi'_\lambda(U_\lambda)=0,\quad \Phi_\lambda(U_\lambda)\le c_\lambda.
$$
\textbf{Proof.} Since
\(\{U_j\}\subset X\) is bounded, we may pass to a subsequence such that
$$
U_j\rightharpoonup U_\lambda \ \mbox{in} \ X, \quad U_{j}\to U_{\lambda} \ \mbox{a.e.} \ x\in\Omega,\qquad j\to\infty,
\eqno(2.15)
$$
where
$$
U_j=(u_j^1,\cdots,u_j^k),\quad U_\lambda=(u_\lambda^1,\cdots,u_\lambda^k).
$$

{\bf Part 1.} We prove \(U_\lambda\neq\mathbf{0}\).
From \(\Phi_\lambda(t l_0 U_0)\mid_{t=1}=\Phi_\lambda(\overline{U})\) (in the proof of Lemma 2.3) and
\(\Phi_\lambda(t l_0 U_0)\mid_{t=0}=\Phi_\lambda(\mathbf{0})\),
we know the map: \([0, 1] \ni t \rightarrow t l_0 U_0\) belongs to \(\Gamma\) in (2.12).
Thus, by \((S_1)\), we have \(c_\lambda\) (in (2.13)) satisfies
\[
\begin{array}{rcl}
    c_\lambda&=&\inf_{\gamma\in\Gamma}\sup_{t\in[0,1]}\Phi_\lambda(\gamma(t)) \\
    \\
    &\le& \sup_{t\in[0,1]} \Phi_\lambda(tl_0 U_0)\\
    \\
    &\le& \sup_{\zeta \ge 0} \Phi(\zeta U_0)\\
    \\
    & <& 2^{-\frac{N-2s}{2s}}\frac{s}{N} S_K^{N/(2s)}
    \le \lambda^{-\frac{N-2s}{2s}}\frac{s}{N} S_K^{N/(2s)}, \quad \forall\lambda\in[1,2].
\end{array}
\]
Thus, we have
\begin{equation}
    \tag{2.16} c_\lambda < \lambda^{-\frac{N-2s}{2s}}\frac{s}{N} S_K^{N/(2s)}, \quad \forall\lambda\in[1,2].
\end{equation}

If \(U_\lambda=\mathbf{0}\), then it follows from \((\Gamma_{1})\), \(\delta\in(1,2)\), (2.15) and Lemma 2.1 that
$$
0\le\left|\int_\Omega \gamma(x)|U_j|^{\delta} dx\right|\le \|\gamma\|_{L^\infty}\|U_j\|_{\delta}^{\delta}\to 0,\qquad j\to\infty.
\eqno(2.17)
$$
By \(U_\lambda=\mathbf{0}\), (2.5), (2.6), (2.15), Lemma 2.1 and \(\|U_j\|_{2^{*}_s}\le\tau_{2^{*}_s}\|U_j\|<\infty\), we have
$$\begin{array}{rcl}
\left|\int_{\Omega}\widetilde{F}(x,U_j)dx\right|&=&\left|\int_{\Omega}\left(\frac{1}{2}\nabla F(x,U_j)\cdot U_j-F(x,U_j)\right)dx\right|\\
\\
&\le& 2^{*}_s \varepsilon \|U_j\|_{2^{*}_s}^{2^{*}_s} + C_\varepsilon\|U_j\|_1+\varepsilon \|U_j\|_{2^{*}_s}^{2^{*}_s} + C_\varepsilon\|U_j\|_1\\
\\
&\to&0,\qquad j\to\infty.
\end{array}
\eqno(2.18)$$
Thus, by (2.14), (2.15), (2.17), (2.18) and \(U_\lambda=\mathbf{0}\), we have
$$\begin{array}{cl}
&\left( \frac{1}{2} - \frac{1}{2^*_s}\right) \lambda\|U_j\|^{2^*_s}_{2^*_s}\\
\\
&=\Phi_\lambda(U_j)
-\frac{1}{2}\langle\Phi'_\lambda(U_j),U_j\rangle-\left(\frac{1}{\delta}-\frac{1}{2}\right)\int_\Omega \gamma(x)|U_j|^{\delta} dx-\lambda\int_{\Omega}\widetilde{F}(x,U_j)dx\\
\\
&\to\lim_{j\to\infty}\Phi_\lambda(U_j)+0+0+0\\
\\
&=c_\lambda>0.
\end{array}
\eqno(2.19)
$$
So we have
$$
\left( \frac{1}{2} - \frac{1}{2^*_s}\right) \lambda\|U_j\|^{2^*_s}_{2^*_s}=\frac{s}{N}\lambda\|U_j\|^{2^*_s}_{2^*_s}\to c_\lambda>0.
\eqno(2.20)
$$
Let
$$
\lim_{j\to\infty}\lambda\|U_j\|^{2^*_s}_{2^*_s}=D\in(0,\infty).
\eqno(2.21)
$$
Then (2.20) implies
$$
\frac{s}{N}D= c_\lambda.
\eqno(2.22)
$$

By \(U_\lambda=\mathbf{0}\), (2.15) and Lemma 2.1, we have
$$
\left|\sum_{i=1}^{k}\mu_i\|u^i_j\|^2_{L^2}\right|\le \left(\max_{i=1,\cdots,k}|\mu_i|\right)\|U_j\|_{2}^{2}\to 0,\qquad j\to\infty.
\eqno(2.23)
$$
By \(U_\lambda=\mathbf{0}\), (2.5), (2.15), Lemma 2.1 and \(\|U_j\|_{2^{*}_s}<\infty\), we have
$$\begin{array}{rcl}
\left|\int_{\Omega}\nabla F(x,U_j)\cdot U_jdx\right|&\le&\int_{\Omega}\left|\nabla F(x,U_j)\cdot U_j\right|dx\\
\\
&\le& 2^{*}_s \varepsilon \|U_j\|_{2^{*}_s}^{2^{*}_s} + C_\varepsilon\|U_j\|_1\\
\\
&\to&0,\qquad j\to\infty.
\end{array}
\eqno(2.24)$$
Hence, by \(U_\lambda=\mathbf{0}\), (2.14), (2.15), Lemma 2.1, (2.17), (2.23) and (2.24), we have
$$\begin{array}{rcl}
    0 \leftarrow \langle\Phi'_\lambda(U_j), U_j \rangle
    &=&\| U_j \|^2- \mathlarger{\sum_{i=1}^{k}}\mu_i\|u^i_j\|^2_{L^2}+\int_\Omega \gamma(x)|U_j|^{\delta} dx \\
    \\
    &&-\lambda\|U_j\|^{2^*_s}_{2^*_s}-\lambda\int_{\Omega}\nabla F(x,U_{j})\cdot U_{j}dx\\
    \\
    &\to&\| U_j \|^2-\lambda\|U_j\|^{2^*_s}_{2^*_s},\qquad j\to\infty.
\end{array}
\eqno(2.25)
$$
Thus, we have
$$
\| U_j \|^2-\lambda\|U_j\|^{2^*_s}_{2^*_s}\rightarrow 0,\qquad  j \rightarrow +\infty.
\eqno(2.26)
$$
It follows from (2.21) and \(\| U_j \|\) is bounded in \(\mathbb{R}\) that
$$
\lim_{j\to\infty}\| U_j \|^2=D.
\eqno(2.27)
$$

By \((1.9)\), we have
$$
    \| U_j \|^2 \ge S_K \| U_j \|_{2^*_s}^2. \eqno(2.28)
$$
Thus, by $\lim_{j\to\infty}\|U_j\|^{2^*_s}_{2^*_s}=\frac{1}{\lambda}D$ in (2.21), (2.27) and (2.28), we get
$$
D \ge \lambda^{-2/2^*_s}S_K D^{2/2^*_s},
$$
it follows from \(D>0\) in (2.21) that
$$
D \ge \lambda^{-\frac{N-2s}{2s}}S_K^{2^*_s/(2^*_s-2)}
= \lambda^{-\frac{N-2s}{2s}}S_K^{N/(2s)},\eqno(2.29)
$$
Hence, (2.29) and (2.22) imply
$$
    c_\lambda=\frac{s}{N}D\ge \lambda^{-\frac{N-2s}{2s}} \frac{s}{N} S_K^{N/(2s)}, \quad \forall\lambda\in[1,2].
    \eqno(2.30)
$$

Now, we get a contradiction between (2.16) and (2.30). Therefore, we have
$$
U_\lambda\neq\mathbf{0}.
$$

{\bf Part 2.} We prove \(\Phi'_\lambda(U_\lambda)=0\).
Since \(\{U_j\}\subset X\) is bounded, Lemma 2.1 yields
\[
\{U_j\}\ \text{is bounded in}\ [L^t(\Omega)]^k,\quad
\forall t\in(1,2^*_s],
\]
and hence
\[
\{|U_j|^{t-2}U_j\}\ \text{is bounded in}\
[L^{t/(t-1)}(\Omega)]^k,\quad \forall t\in(1,2^*_s].
\]

By (2.15), \(U_j\rightharpoonup U_\lambda\) in \(X\) and
\(U_j\to U_\lambda\) a.e. in \(\Omega\). Since \([L^t(\Omega)]^k\) is
reflexive and \(\{U_j\}\) is bounded in \([L^t(\Omega)]^k\), up to a
subsequence (still denoted by \(\{U_j\}\)) we have
\[
U_j\rightharpoonup U_\lambda\quad\text{in}\ [L^t(\Omega)]^k,\quad
\forall t\in(1,2^*_s],
\eqno(2.31)
\]
where the identification of the weak limit with \(U_\lambda\) follows
from the a.e. convergence. Combining this with the pointwise convergence
and the boundedness of \(\{|U_j|^{t-2}U_j\}\) in the reflexive space
\([L^{t/(t-1)}(\Omega)]^k\), we obtain
\[
|U_j|^{t-2}U_j\rightharpoonup|U_\lambda|^{t-2}U_\lambda\quad\text{in}\
[L^{t/(t-1)}(\Omega)]^k,\quad \forall t\in(1,2^*_s],
\eqno(2.32)
\]
as \(j\to\infty\).

Now, (2.31) and Lemma 2.1 implies that
$$
(u^i_j,\varphi^i)_{L^2}\to(u^i,\varphi^i)_{L^2},\quad \forall i=1,\cdots,k,\quad \forall \varphi=(\varphi^1,\cdots,\varphi^k)\in X,\qquad  j \rightarrow +\infty.
\eqno(2.33)
$$
Since \([L^{\delta}(\Omega)]^k\) and \([L^{2^*_s}(\Omega)]^k\) are the duals of \([L^{\delta/(\delta-1)}(\Omega)]^k\) and \([L^{2^*_s/(2^*_s-1)}(\Omega)]^k\) respectively, which together with (2.32) (taking \(t=\delta\) and \(2^*_s\), respectively) and Lemma 2.1 implies that
$$
\int_{\Omega}\gamma(x)|U_j|^{\delta-2}U_j\cdot\varphi \ dx\to\int_{\Omega}\gamma(x)|U_\lambda|^{\delta-2}U_\lambda\cdot\varphi \ dx,\quad \forall \varphi\in X,\qquad  j \rightarrow +\infty,
\eqno(2.34)
$$
and
$$
\int_{\Omega}|U_j|^{2^*_s-2}U_j\cdot\varphi \ dx\to\int_{\Omega}|U_\lambda|^{2^*_s-2}U_\lambda\cdot\varphi \ dx,\quad \forall \varphi\in X,\qquad  j \rightarrow +\infty.
\eqno(2.35)
$$
By (2.5) and the Minkowski inequality
$$(a + b)^l \le 2^{l-1} (a^l + b^l), \quad \forall a,b>0, \ \ \forall l>1,
$$
we have
$$\begin{array}{cl}
|\nabla F(x,U_{j})|^{2^*_s/(2^*_s-1)}&\le \left(2^{*}_s \varepsilon |U_{j}|^{2^{*}_s-1} + C_\varepsilon\right)^{2^*_s/(2^*_s-1)}\\
\\
&\le 2^{2^*_s/(2^*_s-1)-1}\left((2^{*}_s \varepsilon)^{2^*_s/(2^*_s-1)} |U_{j}|^{2^{*}_s} + C_\varepsilon^{2^*_s/(2^*_s-1)}\right).
\end{array}
\eqno(2.36)
$$
Thus, by (2.36) and the facts that \(\{U_j\}\) is bounded in \([L^{2^*_s}(\Omega)]^k\), and \(\Omega\) is bounded, we get
$$
\{\nabla F(x,U_{j})\} \text{ is bounded in } [L^{2^*_s/(2^*_s-1)}(\Omega)]^k.
\eqno(2.37)
$$
Besides, by (2.15) and \(F\in C^1\), we get
$$
\nabla F(x,U_{j}) \to \nabla F(x,U_{\lambda})  \quad \text{a.e. in } \Omega,\qquad  j \rightarrow +\infty.
\eqno(2.38)
$$
Note that \([L^{2^*_s/(2^*_s-1)}(\Omega)]^k\) is a reflexive space, hence (2.37) and (2.38) imply that
\[
\nabla F(x,U_{j})\rightharpoonup \nabla F(x,U_{\lambda})  \quad \text{in } [L^{2^*_s/(2^*_s-1)}(\Omega)]^k,\qquad  j \rightarrow +\infty.
\eqno(2.39)
\]
By (2.39), Lemma 2.1 and the fact that \([L^{2^*_s}(\Omega)]^k\) is the dual of \([L^{2^*_s/(2^*_s-1)}(\Omega)]^k\), we get
$$
    \int_{\Omega} \nabla F(x,U_{j})\cdot\varphi \, dx \rightarrow \int_{\Omega} \nabla F(x,U_{\lambda})\cdot\varphi \, dx,\quad \forall \varphi\in X,\qquad  j \rightarrow +\infty.
    \eqno(2.40)
$$

Thus, by (2.14), (2.15), (2.33)--(2.35), (2.40) and \((\Gamma_1)\), we have that
$$\begin{array}{rcl}
0&\leftarrow&\langle \Phi'_{\lambda}(U_{j}),\varphi\rangle\\
\\
&=&(U_{j},\varphi)- \mathlarger{\sum_{i=1}^{k}}\mu_i(u^i_j,\varphi^i)_{L^2}
+\int_\Omega \gamma(x) |U_j|^{\delta-2}U_j\cdot\varphi  dx\\
\\
&&-\lambda\int_{\Omega}|U_j|^{2^*_s-2}U_j\cdot\varphi dx -\lambda\int_{\Omega} \nabla F(x,U_{j})\cdot\varphi dx \\
\\
&\to& (U_{\lambda},\varphi)- \mathlarger{\sum_{i=1}^{k}}\mu_i(u^i_{\lambda},\varphi^i)_{L^2}
+ \int_\Omega \gamma(x) |U_{\lambda}|^{\delta-2}U_{\lambda}\cdot\varphi  dx\\
\\
&&-\lambda\int_{\Omega}|U_{\lambda}|^{2^*_s-2}U_{\lambda}\cdot\varphi dx -\lambda\int_{\Omega} \nabla F(x,U_{\lambda})\cdot\varphi dx \\
\\
&=&\langle \Phi'_{\lambda}(U_{\lambda}),\varphi\rangle, \quad \forall \varphi=(\varphi^1,\cdots,\varphi^k)\in X, \qquad   \ j \to \infty.
\end{array}
\eqno(2.41)
$$
Therefore, we have
\[
\Phi'_\lambda(U_\lambda)=0.
\tag{2.42}
\]

{\bf Part 3.} We prove \(\Phi_\lambda(U_\lambda)\le c_\lambda\).
By (2.15), \(U_j\to U_\lambda\) a.e. in \(\Omega\), and by \((F_1)\),
\(\widetilde F(x,U_j)\to\widetilde F(x,U_\lambda)\) a.e. in \(\Omega\).
Since \(\widetilde F(x,U_j)\ge0\) a.e. in \(\Omega\) for every \(j\) by
\((F_4)\), the Fatou lemma gives
\[
\int_\Omega\widetilde F(x,U_\lambda)\,dx
\le\liminf_{j\to\infty}\int_\Omega\widetilde F(x,U_j)\,dx.
\tag{2.43}
\]
By (2.31) and the compact embedding
\(X\hookrightarrow[L^\delta(\Omega)]^k\), \(\delta\in(1,2)\), we have
\[
\int_\Omega\gamma|U_j|^\delta\,dx\to
\int_\Omega\gamma|U_\lambda|^\delta\,dx.
\]
By the weak lower semicontinuity of the \(L^{2^*_s}\) norm,
\[
\|U_\lambda\|_{2^*_s}^{2^*_s}
\le\liminf_{j\to\infty}\|U_j\|_{2^*_s}^{2^*_s}.
\]
Combining these with (2.14),
\(2^*_s>2\), (2.43) and (2.42), we obtain
$$\begin{array}{rcl}
c_\lambda&=&\lim_{j\to\infty}\left(\Phi_\lambda(U_j)
-\frac{1}{2}\langle\Phi'_\lambda(U_j),U_j\rangle\right)\\
\\
&=&\lim_{j\to\infty}\left(\left(\frac{1}{\delta}-\frac{1}{2}\right)\int_\Omega \gamma(x) |U_j|^{\delta}dx
+\lambda\left(\frac{1}{2}-\frac{1}{2^*_s}\right)\|U_j\|_{2^*_s}^{2^*_s}
+\lambda\int_{\Omega}\widetilde{F}(x,U_{j})dx\right)\\
\\
&\ge&\left(\frac{1}{\delta}-\frac{1}{2}\right)\int_\Omega \gamma(x) |U_{\lambda}|^{\delta}dx
+\lambda\left(\frac{1}{2}-\frac{1}{2^*_s}\right)\|U_{\lambda}\|_{2^*_s}^{2^*_s}
+\lambda\int_{\Omega}\widetilde{F}(x,U_{\lambda})dx\\
\\
&=&\Phi_\lambda(U_\lambda)-\frac{1}{2}\langle\Phi'_\lambda(U_\lambda),U_\lambda\rangle\\
\\
&=&\Phi_\lambda(U_\lambda).
\end{array}
\eqno(2.44)
$$
Therefore, the proof is finished.$\quad\Box$\\

Now, Lemmas 2.4 and 2.5 imply that there is \(\{(\lambda_j,U_j)\}\subset[1,2]\times X\) such that
$$
\{\lambda_j\} \ \mbox{is decreasing}, \quad
\lambda_j\to1,\quad U_j\neq0,\quad
\Phi_{\lambda_j}(U_j)\le c_{\lambda_j},\quad \Phi'_{\lambda_j}(U_j)=0.
\eqno(2.45)
$$

In what follows, we establish the boundedness of PS sequence \(\{U_j\}\) for the approximating functionals \(\Phi_{\lambda_j}\). Namely, the following lemma.

\medskip
\noindent\textbf{Lemma 2.6 (Boundedness of the PS sequence).} {\it Under the conditions of Theorem $1.1$, the sequence \(\{U_j\}\subset X\) obtained in \((2.45)\) is bounded.}\\
\textbf{Proof.} Arguing by contradiction, suppose that \(\|U_j\|\rightarrow\infty\).
Let
$$
(v_j^1,\cdots,v_j^k)=V_j:=\frac{U_j}{\|U_j\|}=\left(\frac{u_j^1}{\|U_j\|},\cdots,\frac{u_j^k}{\|U_j\|}\right),
\eqno(2.46)
$$
then
$$\|V_j\|=1.$$
Thus, up to a subsequence (still denoted by \(\{V_j\}\)), we have
$$
V_j\rightharpoonup V \ \mbox{ in} \ X, \quad
V_j\to V \ \mbox{ a.e.} \ x\in\Omega, \qquad   \ j \to \infty.
\eqno(2.47)
$$

Since \(\Phi'_{\lambda_j}(U_j)=0\), \(\Phi_{\lambda_j}(U_j)\le c_{\lambda_j}\)
by (2.45), and \(\{c_{\lambda_j}\}\) is uniformly bounded above by (2.16)
and \(\lambda_j\in[1,2]\), there exists \(\widetilde{\kappa}>0\) such that
\[
\Phi_{\lambda_j}(U_j) - \frac{1}{2} \langle \Phi_{\lambda_j}'(U_j), U_j \rangle
\le \widetilde{\kappa} (1 + \| U_j \|).
\tag{2.48}
\]
By \((F_4)\) and \(\gamma(x)\ge0\), \(1<\delta<2\), \(\lambda_j\in[1,2]\),
we obtain
$$
\begin{aligned}
&\Phi_{\lambda_j}(U_j) - \frac{1}{2} \langle \Phi_{\lambda_j}'(U_j), U_j \rangle\\
\\
&=\left(\frac{1}{\delta}-\frac{1}{2}\right)\int_\Omega \gamma(x) |U_j|^{\delta}dx
+\left( \frac{1}{2} - \frac{1}{2^*_s}\right)\lambda_j \| U_j \|_{2^*_s}^{2^*_s}
+\lambda_j\int_{\Omega}\widetilde{F}(x,U_j)dx\\
\\
&\geq \frac{s}{N} \| U_j \|_{2^*_s}^{2^*_s}.
\end{aligned}
\eqno(2.49)
$$
Combining (2.48) and (2.49) yields
\[
\|U_j\|_{2^*_s}^{2^*_s}
\le \frac{N}{s}\widetilde{\kappa}(1+\|U_j\|),
\eqno(2.50)
\]

By \((2.5)\) and Lemma 2.1, we have
$$
\int_{\Omega}\nabla F(x,U_{j})\cdot U_{j}dx\leq 2^{*}_s \varepsilon \|U_j\|_{2^{*}_s}^{2^{*}_s} + C_\varepsilon\tau_1\|U_j\|.
\eqno(2.51)
$$
Thus, the facts \(\Phi'_{\lambda_j}(U_j)=0\), \(\gamma(x)\ge0\), (2.50) and (2.51) imply there exists \(\widetilde{\kappa_0}>0\) such that
$$\begin{array}{rcl}
\|U_j\|^2&=&\mathlarger{\sum_{i=1}^{k}}\mu_i\|u^i_j\|_{L^2}^2-\int_\Omega \gamma(x) |U_j|^{\delta}dx+\lambda_j \|U_j\|_{2^*_s}^{2^*_s}+\lambda_j\int_{\Omega}\nabla F(x,U_{j})\cdot U_{j}dx\\
\\
&\le&\left(\mathlarger{\max_{i=1,\cdots,k}}|\mu_i|\right)\|U_j\|^{2}_2+\widetilde{\kappa_0}(1 + \| U_j \| ).
\end{array}
\eqno(2.52)
$$

We claim that neither \(V=0\) nor \(V\neq0\) is possible; this
contradiction will show that \(\{U_j\}\) is bounded.

{\bf Case 1.} If \(V=\mathbf{0}\).
By (2.47), (2.52), \(V=\mathbf{0}\), Lemma 2.1 and \(\|U_j\|\rightarrow\infty\), we get
$$
1\le\left(\max_{i=1,\cdots,k}|\mu_i|\right)\|V_j\|^{2}_2
+\frac{\widetilde{\kappa_0}(1 + \| U_j\|)}{\|U_j\|^2}\to0, \qquad   \ j \to \infty.
$$
Therefore, \(V=\mathbf{0}\) is impossible.

{\bf Case 2.} If \(V\neq\mathbf{0}\).
By \(\Phi'_{\lambda_j}(U_j)=0\) and (2.5) (taking \(\varepsilon>0\) sufficiently small), we have
$$\begin{array}{cl}
\|U_{j}\|^2&=\mathlarger{\sum_{i=1}^{k}}\mu_i\|u^i_j\|^2_{L^2}-\int_\Omega \gamma(x)|U_{j}|^{\delta} dx
+\lambda_j\|U_{j}\|_{2^*_s}^{2^*_s}+\lambda_j\int_{\Omega}\nabla F(x,U_{j})\cdot U_{j} dx\\
\\
&\ge-\mathlarger{\max_{i=1,\cdots,k}}|\mu_i|\tau^2_2\|U_{j}\|^2 -\|\gamma\|_{L^\infty}\tau_\delta^\delta\|U_{j}\|^\delta\\
\\
&\quad +\|U_{j}\|_{2^*_s}^{2^*_s}
-\varepsilon(2\times2^{*}_s)\|U_{j}\|_{2^*_s}^{2^*_s}-2C_\varepsilon\tau_1\|U_{j}\|.\\
\\
&\ge-\mathlarger{\max_{i=1,\cdots,k}}|\mu_i|\tau^2_2\|U_{j}\|^2-\|\gamma\|_{L^\infty}\tau_\delta^\delta\|U_{j}\|^\delta
+\frac{1}{2}\|U_{j}\|_{2^*_s}^{2^*_s}-2C_\varepsilon\tau_1\|U_{j}\|.
\end{array}
\eqno(2.53)
$$
By (2.53), \(V\neq\mathbf{0}\) and \(V_j=\frac{U_j}{\|U_j\|}\),
$$
1\ge
-\max_{i=1,\cdots,k}|\mu_i|\tau^2_2-\frac{\|\gamma\|_{L^\infty}\tau_\delta^\delta}{\|U_{j}\|^{2-\delta}}
+\frac{1}{2}\int_{\Omega}|U_j|^{2^*_s-2}|V_{j}|^2dx-\frac{2C_\varepsilon\tau_1}{\|U_{j}\|}.
$$
It follows from \(\|U_j\|\to\infty\) that
$$
1\ge
-\max_{i=1,\cdots,k}|\mu_i|\tau^2_2
+\lim_{j\to\infty}\frac{1}{2}\int_{\Omega}|U_j|^{2^*_s-2}|V_{j}|^2dx.
\eqno(2.54)
$$
Since \(V\neq\mathbf{0}\) and \(\|U_j\|\to\infty\), we have
\[
|U_j|=|V_j|\cdot\|U_j\|\to\infty.
\]
Moreover,
\[
\int_\Omega|U_j|^{2^*_s-2}|V_j|^2\,dx
=\|U_j\|^{2^*_s-2}\int_\Omega|V_j|^{2^*_s}\,dx.
\eqno(2.55)
\]
Since \(V_j\rightharpoonup V\neq0\) in \(X\) and the embedding
\(X\hookrightarrow L^{2^*_s}(\Omega)\) is continuous, we have
\(V_j\rightharpoonup V\) in \(L^{2^*_s}(\Omega)\). By the weak lower
semicontinuity of the \(L^{2^*_s}\) norm,
\[
\liminf_{j\to\infty}\int_\Omega|V_j|^{2^*_s}\,dx
\ge\int_\Omega|V|^{2^*_s}\,dx>0.
\]
Hence, since \(\|U_j\|\to\infty\),
\[
\int_\Omega|U_j|^{2^*_s-2}|V_j|^2\,dx
=\|U_j\|^{2^*_s-2}\int_\Omega|V_j|^{2^*_s}\,dx\to+\infty.
\eqno(2.56)
\]
Now, by (2.54) and (2.56), we obtain a contradiction.
Thus, \(V\neq\mathbf{0}\) is impossible.

Therefore, \(\{U_j\}\subset X\) is bounded.$\quad\Box$

\section*{3. Proofs of Main Results}

\noindent\textbf{Proof of Theorem 1.1.} From the boundedness of \(\{U_j\}\subset X\) in Lemma 2.6, up to a subsequence (still denoted by \(\{U_j\}\)), we have
$$
U_j\rightharpoonup U \ \mbox{in} \ X,\quad U_{j}\to U, \ a.e. \ x\in\Omega,\qquad j\to\infty.
\eqno(3.1)
$$
Next, we divide the proof into three parts.

{\bf Part 1.} We prove \(\Phi'(U)=0\).
Similar to (2.35) and (2.40), by (3.1) and Lemma 2.1, we have
$$
\int_{\Omega}|U_j|^{2^*_s-2}U_j\cdot\varphi \ dx\to\int_{\Omega}|U|^{2^*_s-2}U \cdot\varphi \ dx,\quad \forall \varphi\in X,\qquad  j \rightarrow +\infty,
\eqno(3.2)
$$
and
$$
\int_{\Omega} \nabla F(x,U_{j})\cdot\varphi \, dx \rightarrow \int_{\Omega} \nabla F(x,U)\cdot\varphi \, dx,\quad \forall \varphi\in X,\qquad  j \rightarrow +\infty.
\eqno(3.3)
$$
Thus, by $\Phi'_{\lambda_j}(U_j)=0$ in (2.45), (3.2), (3.3) and \(\lambda_{j}\to1\), we have
$$\begin{array}{cl}
&\langle\Phi'(U_{j}),\varphi\rangle\\
\\
&=
\langle\Phi'_{\lambda_{j}}(U_{j}),\varphi\rangle
+(\lambda_{j}-1)\int_{\Omega}|U_j|^{2^*_s-2} U_j\cdot \varphi dx+(\lambda_{j}-1)\int_{\Omega}\nabla F(x,U_{j}) \cdot \varphi dx\\
\\
&\to0,\quad\forall \varphi\in
X.
\end{array}
\eqno(3.4)
$$
Similar to (2.41) in Part 2 of Lemma 2.5, by (3.1) and (3.4), we have
$$
0=\lim_{j\to \infty}\langle\Phi'(U_j),\varphi\rangle=\langle\Phi'(U),\varphi\rangle, \quad \forall \varphi\in X.
$$
Therefore, we have
$$\Phi'(U)=0.$$

{\bf Part 2.} We prove \(U\) is a nontrivial solution of the system (1.2)
By (2.45), (2.7) and Lemma 2.1, we have
$$
\begin{array}{cl}
&\|U_j\|^2-\mathlarger{\sum_{i=1}^{k}}\mu_i\|u^i_j\|_{L^2}^2
+\int_\Omega \gamma(x) |U_j|^{\delta}dx\\
\\
&=\lambda_j\| U_j \|_{2^*_s}^{2^*_s}+\lambda_j\int_{\Omega}\nabla F(x,U_{j}) \cdot U_j dx\\
\\
&\le2\| U_j \|_{2^*_s}^{2^*_s}+ 4\varepsilon\|U_j\|_2^2+ 2\times2^{*}_s\delta_\varepsilon\|U_j\|_{2^{*}_s}^{2^{*}_s} \\
\\
&\le  \varepsilon4\tau_2^2\|U_j\|^2+\left(2+2\times2^{*}_s\delta_\varepsilon\right)\| U_j \|_{2^*_s}^{2^*_s-2}\tau_{2^*_s}^{2}\|U_j\|^2.
\end{array}
$$
On the other hand, \(\gamma(x)\ge0\) and Lemma 2.2 imply
$$
\|U_j\|^2-\sum_{i=1}^{k}\mu_i\|u_i\|_{L^2}^2+\int_\Omega \gamma(x) |U_j|^{\delta}dx
\ge\overline{D}\|U_j\|^2.
$$
Hence, we have
$$
\overline{D}\|U_j\|^2\le \varepsilon4\tau_2^2\|U_j\|^2
+\left(2+2\times2^{*}_s\delta_\varepsilon\right)\| U_j \|_{2^*_s}^{2^*_s-2}\tau_{2^*_s}^{2}\|U_j\|^2, \quad \forall\varepsilon>0.
$$
Since \(U_j\neq0\), dividing by \(\|U_j\|^2\) gives
\[
\overline{D}\le \varepsilon4\tau_2^2
+\bigl(2+2\cdot 2^*_s\delta_\varepsilon\bigr)
\tau_{2^*_s}^2\|U_j\|_{2^*_s}^{2^*_s-2}.
\]
Choosing \(\varepsilon>0\) small enough so that
\(\varepsilon4\tau_2^2<\overline{D}/2\), we obtain
\[
\|U_j\|_{2^*_s}\ge\widetilde{\omega}
:=\left(
\frac{\overline{D}/2}
{\bigl(2+2\cdot 2^*_s\delta_\varepsilon\bigr)\tau_{2^*_s}^2}
\right)^{1/(2^*_s-2)}>0,
\tag{3.5}
\]
where \(\widetilde{\omega}\) is independent of \(j\).

Note that \(\|U_j\|<\infty\) and Lemma 2.1 imply that \(\|U_j\|_{2^{*}_s}\le\tau_{2^{*}_s}\|U_j\|<\infty\). Thus, by (3.5) and \(\lambda_j\to1\), we can let
\begin{equation}
\lim_{j\to\infty}\|U_j\|_{2^*_s}^{2^*_s}
=\lim_{j\to\infty}\lambda_j\|U_j\|_{2^*_s}^{2^*_s}
=\widetilde{D}\in[\widetilde{\omega}^{2^*_s},\infty).
\tag{3.6}
\end{equation}

Assume, for contradiction, that \(U=\mathbf{0}\). Then
\(U_j\rightharpoonup0\) in \(X\) and \(U_j\to0\) a.e. in \(\Omega\).
Similar to Part 1 in the proof of Lemma 2.5 (see (2.19) and (2.26)), we also have
\begin{equation}
\begin{aligned}
&\Bigl(\frac12-\frac1{2^*_s}\Bigr)\lambda_j\|U_j\|_{2^*_s}^{2^*_s}\\
&=\Phi_{\lambda_j}(U_j)-\tfrac12\langle\Phi'_{\lambda_j}(U_j),U_j\rangle
-\Bigl(\frac1\delta-\frac12\Bigr)\int_\Omega\gamma|U_j|^\delta\,dx
-\lambda_j\int_\Omega\widetilde F(x,U_j)\,dx\\
&\to c_{\lambda_j}+0+0+0,
\end{aligned}
\tag{3.7}
\end{equation}
and
$$
\| U_j \|^2-\lambda_j\|U_j\|^{2^*_s}_{2^*_s}\rightarrow 0,\qquad  j \rightarrow +\infty.
\eqno(3.8)
$$
By \(\lambda_j\to1\), (3.6) and (3.7),
\begin{equation}
\frac{s}{N}\widetilde{D}= c_{\lambda_j}.
\tag{3.9}
\end{equation}
By (3.6) and (3.8),
\[
\|U_j\|^2\to\widetilde{D}.
\]
Combining this with (3.6) and the Sobolev inequality
\(\|U_j\|^2\ge S_K\|U_j\|_{2^*_s}^2\) in (1.9),
we obtain
\[
\widetilde{D}\ge S_K\widetilde{D}^{2/2^*_s},
\tag{3.10}
\]
hence
\[
\widetilde{D}\ge S_K^{N/(2s)}.
\tag{3.11}
\]
Combining (3.9), (3.11) and $\lambda_j\in[1,2]$, we have
\begin{equation}
c_{\lambda_j}=\frac{s}{N}\widetilde{D}
\ge\frac{s}{N}S_K^{N/(2s)}.
\tag{3.12}
\end{equation}
On the other hand, as in Part 1 of the proof of Lemma 2.5 (see (2.16)), we have
\begin{equation}
c_{\lambda_j}< {\lambda_j}^{-\frac{N-2s}{2s}}\frac{s}{N} S_K^{N/(2s)}\le \frac{s}{N} S_K^{N/(2s)}, \quad \forall\lambda_j\in[1,2].
\tag{3.13}
\end{equation}

Now (3.12) and (3.13) contradict each other. Hence
$$
U\neq\mathbf{0}.
$$
Therefore, from Part 1, \(U\) is a nontrivial solution of the system (1.2).

{\bf Part 3.} We prove that the system (1.2) has a ground state solution.
Let
\[
\widetilde{K}:=\{U\in X:\ \Phi'(U)=0\},
\qquad
\widetilde{c}:=\inf\left\{\Phi(U):\ U\in \widetilde{K}\setminus\{\mathbf{0}\}\right\}.
\]
Let \(\{U_j\}\subset\widetilde{K}\setminus\{\mathbf{0}\}\) be a
minimizing sequence with
\begin{equation}
\Phi'(U_j)=0,\quad \Phi(U_j)\to\widetilde{c}.
\tag{3.14}
\end{equation}
Thus, by the same argument as in Lemma 2.6 (with \(\lambda_j\) replaced by
\(1\)), \(\{U_j\}\) is bounded in \(X\). Up to a subsequence, there
exists \(\widetilde U\in X\) with
\[
U_j\rightharpoonup\widetilde U \ \mbox{in} \ X,\quad
U_j\to\widetilde U \ \mbox{a.e. in} \  \Omega.
\]
Arguing as in Part 1 and
Part 2 above, we obtain
\[
\widetilde U\neq0,\quad
\Phi'(\widetilde U)=0.
\]

Since \(\widetilde F(x,U_j)\ge0\) by \((F_4)\) and
\(\widetilde F(x,U_j)\to\widetilde F(x,\widetilde{U})\) a.e. in
\(\Omega\) by (3.14) and \((F_1)\), it follows from the Fatou lemma, the weak lower
semicontinuity of \(\|\cdot\|_{2^*_s}\), the compact embedding
\(X\hookrightarrow[L^\delta(\Omega)]^k\) and \(\Phi'(\widetilde{U})=0\)
that
$$\begin{array}{rcl}
\widetilde{c}&=&\lim_{j\rightarrow\infty}\left(\Phi(U_{j})-\frac{1}{2}\langle\Phi'(U_{j}),U_{j}\rangle\right)\\
\\
&=&\lim_{j\to\infty}\left(\left(\frac{1}{\delta}-\frac{1}{2}\right)\int_\Omega \gamma(x) |U_j|^{\delta}dx
+\left(\frac{1}{2}-\frac{1}{2^*_s}\right)\|U_j\|_{2^*_s}^{2^*_s}
+\int_{\Omega}\widetilde{F}(x,U_{j})dx\right)\\
\\
&\geq&\left(\frac{1}{\delta}-\frac{1}{2}\right)\int_\Omega \gamma(x) |\widetilde{U}|^{\delta}dx
+\frac{s}{N}\|\widetilde{U}\|_{2^*_s}^{2^*_s}
+\int_{\Omega}\widetilde{F}(x,\widetilde{U})dx\\
\\
&=&\Phi(\widetilde{U})-\frac{1}{2}\langle\Phi'(\widetilde{U}),\widetilde{U}\rangle\\
\\
&=&\Phi(\widetilde{U})\geq \widetilde{c},\end{array}
\eqno(3.15)
$$
which implies \(\Phi(\widetilde{U})=\widetilde{c}\). Since \(\widetilde U\neq0\),
\(\widetilde U\) is a ground state solution.

By (3.15), \(\widetilde F\ge0\), \(\widetilde U\neq0\) and \(\delta\in(1,2)\),
\[
\widetilde{c}=\Phi(\widetilde{U})
=\Bigl(\frac{1}{\delta}-\frac{1}{2}\Bigr)\int_\Omega\gamma|\widetilde{U}|^\delta
+\frac{s}{N}\|\widetilde{U}\|_{2^*_s}^{2^*_s}
+\int_\Omega\widetilde F(x,\widetilde{U})>0.
\]
Therefore, the proof is finished. $\quad\Box$

\medskip
\noindent\textbf{Proof of Theorem 1.2.} Recall that \(u_0\) is the positive first
eigenfunction of \(-\mathcal L_{K_1}\) associated with \(\sigma_1^1\),
normalized by \(\|u_0\|_{L^2}=1\), and \(U_0=(u_0,0,\dots,0)\) (see Notation 1.1).
Then
\[
\|U_0\|^2=\|u_0\|_{\widetilde E_1}^2=\sigma_1^1\|u_0\|_{L^2}^2=\sigma_1^1,
\qquad
\|U_0\|_{L^2}^2=1.
\]
Define
$$
A:=\sigma_1^1-\mu_1>0,
\qquad
B:=\|U_0\|_{2^*_s}^{2^*_s}=\|u_0\|_{L^{2^*_s}}^{2^*_s},
\eqno(3.16)
$$
where \(A>0\) follows from \((M_1)\).

\smallskip

{\bf Step 1 (The pure critical case).} Let \(\gamma\equiv0\), \(F\equiv0\).
Let \(\Phi_0\) denote the functional corresponding to \(\gamma\equiv0\),
\(F\equiv0\).

For \(\zeta\ge0\), by (3.16) and the definition of $\Phi$ in (1.7),
$$
\Phi_0(\zeta U_0)=\frac{A}{2}\zeta^2-\frac{B}{2^*_s}\zeta^{2^*_s}.
\eqno(3.17)
$$
Differentiating gives
\[
\frac{d}{d\zeta}\Phi_0(\zeta U_0)=A\zeta-B\zeta^{2^*_s-1}=0
\iff \zeta^{2^*_s-2}=\frac{A}{B}.
\]
Hence the maximum is attained at
$$
\zeta_0:=\left(\frac{A}{B}\right)^{\frac{1}{2^*_s-2}}
=\left(\frac{A}{B}\right)^{\frac{N-2s}{4s}},
\eqno(3.18)
$$
and
\[
\sup_{\zeta\ge0}\Phi_0(\zeta U_0)
=\left(\frac{A}{B}\right)^{\frac{N-2s}{2s}}
\left(\frac{A}{2}-\frac{A}{2^*_s}\right)
=\frac{s}{N}A\left(\frac{A}{B}\right)^{\frac{N-2s}{2s}}.
\]
Introduce
$$
Q(u_0):=\frac{A}{B^{2/2^*_s}}
=\frac{\sigma_1^1-\mu_1}{\|u_0\|_{L^{2^*_s}}^2}>0.
\eqno(3.19)
$$
Then
$$
\sup_{\zeta\ge0}\Phi_0(\zeta U_0)
=\frac{s}{N}\bigl[Q(u_0)\bigr]^{N/(2s)}.
\eqno(3.20)
$$
By (1.10), we have
\[
Q(u_0)<2^{-\frac{N-2s}{N}}S_K
\quad\Longrightarrow\quad
\bigl[Q(u_0)\bigr]^{N/(2s)}
<2^{-\frac{N-2s}{2s}}S_K^{N/(2s)}.
\]
Hence, it follows from (3.20) that
$$
\sup_{\zeta\ge0}\Phi_0(\zeta U_0)
=\frac{s}{N}\bigl[Q(u_0)\bigr]^{N/(2s)}
<2^{-\frac{N-2s}{2s}}\frac{s}{N}S_K^{N/(2s)}.
\eqno(3.21)
$$
So \((S_1)\) holds in the pure critical case.

\smallskip

{\bf Step 2 (The small-perturbation case).}
Recall that \(\Phi_0\) denotes the functional corresponding to \(\gamma\equiv0\),
\(F\equiv0\).
For \(\zeta\ge0\), using \(\|\gamma\|_{L^\infty}\le\varepsilon_0\) and
\(|\nabla F(x,U)|\le\varepsilon_0|U|+\varepsilon_0|U|^{2^*_s-1}\) in
(1.11), we have
$$
\begin{aligned}
\bigl|\Phi(\zeta U_0)-\Phi_0(\zeta U_0)\bigr|
&\le\frac{1}{\delta}\int_\Omega\gamma|\zeta U_0|^\delta\,dx
+\int_\Omega|F(x,\zeta U_0)|\,dx\\
\\
&\le\frac{\varepsilon_0}{\delta}\|U_0\|_\delta^\delta\zeta^\delta
+\frac{\varepsilon_0}{2}\|U_0\|_{2}^2\zeta^2
+\frac{\varepsilon_0}{2^*_s}\|U_0\|_{2^*_s}^{2^*_s}\zeta^{2^*_s}\\
\\
&\le C(\varepsilon_0)\bigl(\zeta^\delta+\zeta^2+\zeta^{2^*_s}\bigr).
\end{aligned}
\eqno(3.22)
$$
where
\[
C(\varepsilon_0)
:=\max\Bigl\{\frac{\varepsilon_0}{\delta}\|U_0\|_\delta^\delta,\;
\frac{\varepsilon_0}{2}\|U_0\|_{2}^2,\;
\frac{\varepsilon_0}{2^*_s}\|U_0\|_{2^*_s}^{2^*_s}\Bigr\}
\]
depends on \(U_0\) but satisfies
$$
C(\varepsilon_0)\to0^+,\quad
\varepsilon_0\to0^+.
\eqno(3.23)
$$

By (3.21) in Step 1,
there exists \(\eta>0\) such that
$$
\sup_{\zeta\ge0}\Phi_0(\zeta U_0)
\le2^{-\frac{N-2s}{2s}}\frac{s}{N}S_K^{N/(2s)}-\eta.
\eqno(3.24)
$$
By (2.9) (with \(\lambda=1\)),
\[
\Phi(\zeta U_0)\to-\infty,\quad
\zeta\to+\infty.
\]
Since \(\Phi(\mathbf{0})=0\), there exists
\(\zeta_*>0\), depending only on \(U_0\), such that
$$
\sup_{\zeta\ge0}\Phi(\zeta U_0)=\sup_{\zeta\in[0,\zeta_*]}\Phi(\zeta U_0).
\eqno(3.25)
$$
By (3.22),
$$
\sup_{\zeta\in[0,\zeta_*]}\bigl|\Phi(\zeta U_0)-\Phi_0(\zeta U_0)\bigr|
\le C(\varepsilon_0)\bigl(\zeta_*^\delta+\zeta_*^2+\zeta_*^{2^*_s}\bigr).
\eqno(3.26)
$$
Choosing \(\varepsilon_0>0\) small enough so that
\[
C(\varepsilon_0)\bigl(\zeta_*^\delta+\zeta_*^2+\zeta_*^{2^*_s}\bigr)<\eta/2,
\]
which is possible by (3.23), and by  (3.24)--(3.26), we obtain
\[
\begin{aligned}
\sup_{\zeta\ge0}\Phi(\zeta U_0)
&=\sup_{\zeta\in[0,\zeta_*]}\Phi(\zeta U_0)\\
&\le\sup_{\zeta\in[0,\zeta_*]}\Phi_0(\zeta U_0)+\eta/2\\
&\le\sup_{\zeta\ge0}\Phi_0(\zeta U_0)+\eta/2\\
&\le2^{-\frac{N-2s}{2s}}\frac{s}{N}S_K^{N/(2s)}-\eta+\eta/2\\
&<2^{-\frac{N-2s}{2s}}\frac{s}{N}S_K^{N/(2s)}.
\end{aligned}
\]
Therefore \((S_1)\) holds in the small-perturbation case.

\smallskip
In both cases, condition \((S_1)\) holds with \(U_0=(u_0,0,\dots,0)\in X\).
Together with \((K_1)\)--\((K_3)\), \((\Gamma_1)\), \((F_1)\)--\((F_4)\) and
\((M_1)\), all the hypotheses of Theorem 1.1 are satisfied. Hence, by
Theorem 1.1, the system (1.2) admits a ground state solution
\(\widetilde U\in X\setminus\{0\}\) such that
\(\Phi(\widetilde U)=\widetilde c>0.\) \(\Box\)

\medskip
\noindent\textbf{Proof of Theorem 1.3.}
Let \(u_0\) and \(U_0=(u_0,0,\dots,0)\) be as in Notation 1.1, and let
\(A,B,Q(u_0),\zeta_0\) be as in the proof of Theorem 1.2. By Step 1 of
the proof of Theorem 1.2, the maximum of \(\Phi_0\) along the ray
\(\zeta\mapsto\zeta U_0\) is
$$
\sup_{\zeta\ge0}\Phi_0(\zeta U_0) =\frac{s}{N}\bigl[Q(u_0)\bigr]^{N/(2s)},
\eqno(3.27)
$$
attained at \(\zeta_0\).
Define
$$
\eta_0
:=2^{-\frac{N-2s}{2s}}\frac{s}{N}S_K^{N/(2s)}
-\frac{s}{N}\bigl[Q(u_0)\bigr]^{N/(2s)}>0.
\eqno(3.28)
$$
Indeed, by (1.10) and the definitions of \(A\), \(B\), \(Q(u_0)\), we
have \(Q(u_0)<2^{-\frac{N-2s}{N}}S_K\). Raising both sides to the power
\(N/(2s)\) and multiplying by \(s/N\) gives \(\eta_0>0\). By the
definition of \(\eta_0\) and (3.27),
$$
\sup_{\zeta\ge0}\Phi_0(\zeta U_0)
=2^{-\frac{N-2s}{2s}}\frac{s}{N}S_K^{N/(2s)}-\eta_0.
\eqno(3.29)
$$

Since \(F\ge0\) by \((F_1)\), for all \(\zeta\ge0\),
$$
\Phi(\zeta U_0)
\le \Phi_0(\zeta U_0)
+\frac{1}{\delta}\zeta^\delta\int_\Omega\gamma(x)|u_0|^\delta\,dx.
\eqno(3.30)
$$
Define
$$
C_\gamma:=\int_\Omega\gamma(x)|u_0|^\delta\,dx,
\qquad
h(\zeta):=\Phi_0(\zeta U_0)+\frac{C_\gamma}{\delta}\zeta^\delta,
\quad \zeta\ge0.
\eqno(3.31)
$$
Then (3.30) implies
$$
\sup_{\zeta\ge0}\Phi(\zeta U_0)\le \sup_{\zeta\ge0}h(\zeta).
\eqno(3.32)
$$
Since \(\delta<2<2^*_s\), we have \(h(\zeta)\to-\infty\) as
\(\zeta\to+\infty\). Hence \(h\) attains its maximum at some finite
\(\zeta_*\ge0\).

We first show that, for \(C_\gamma\) sufficiently small, the maximizer
\(\zeta_*\) is controlled independently of \(C_\gamma\). The maximizer
satisfies \(h'(\zeta_*)=0\). By (3.17) and (3.31),
$$
h'(\zeta)=A\zeta-B\zeta^{2^*_s-1}+C_\gamma\zeta^{\delta-1},\qquad \zeta\ge0.
\eqno(3.33)
$$
Define
$$
G(\zeta,C_\gamma):=h'(\zeta)=A\zeta-B\zeta^{2^*_s-1}
+C_\gamma\zeta^{\delta-1},\qquad \zeta\ge0.
\eqno(3.34)
$$
Since \(\zeta_0\) is the maximum point of \(\Phi_0(\zeta U_0)\) in the
pure critical case, we have
$$
A\zeta_0-B\zeta_0^{2^*_s-1}=0.
\eqno(3.35)
$$
Hence \(G(\zeta_0,0)=0\) by (3.34) and (3.35). Also, \(G\) is \(C^1\), and
$$
\partial_\zeta G(\zeta_0,0)
=A-B(2^*_s-1)\zeta_0^{2^*_s-2}
=A-A(2^*_s-1)
=-A(2^*_s-2)<0,
$$
where we used \(B\zeta_0^{2^*_s-2}=A\), which follows directly from the
definition of \(\zeta_0\). Thus \(\partial_\zeta G(\zeta_0,0)\neq0\).
By the implicit function theorem, there exist \(\varepsilon_1'>0\) and a
\(C^1\)-map
\[
\zeta_*:[0,\varepsilon_1')\to(0,\infty)
\]
such that
\[
\zeta_*(0)=\zeta_0;\qquad
G(\zeta_*(C_\gamma),C_\gamma)=0,
\quad \forall C_\gamma\in[0,\varepsilon_1').
\]
By continuity, shrinking \(\varepsilon_1'\) if necessary, we may assume
$$
0\le C_\gamma<\varepsilon_1'
\quad\Longrightarrow\quad
0<\zeta_*(C_\gamma)\le 2\zeta_0.
\eqno(3.36)
$$
Moreover, for \(C_\gamma\ge0\) small, this positive zero is the unique
positive zero of \(G(\cdot,C_\gamma)\). Indeed, by (3.34), we can write
$$
G(\zeta,C_\gamma)=\zeta^{\delta-1}f(\zeta),
\qquad
f(\zeta):=A\zeta^{2-\delta}-B\zeta^{2^*_s-\delta}+C_\gamma,\qquad \zeta\ge0.
\eqno(3.37)
$$
Since \(2<\delta<2^*_s\), we have
\[
f'(\zeta)=\zeta^{1-\delta}
\bigl[A(2-\delta)-B(2^*_s-\delta)\zeta^{2^*_s-2}\bigr],\qquad \zeta\ge0.
\]
Thus \(f'(\zeta)=0\) is equivalent to
\[
\zeta^{2^*_s-2}
=\frac{A(2-\delta)}{B(2^*_s-\delta)}>0,
\]
which has exactly one positive solution \(\zeta_f>0\). Consequently, \(f'\) is
positive on \((0,\zeta_f)\) and negative on \((\zeta_f,+\infty)\). Hence \(f\) first increases and then decreases on
\((0,\infty)\). Together with
\[
f(0)=C_\gamma\ge0,\qquad f(\zeta)\to-\infty \ \mbox{as} \ \zeta\to+\infty,
\]
this shows that \(f\)
has exactly one positive zero. Therefore, (3.37) implies \(G(\cdot,C_\gamma)\) has
exactly one positive zero. Consequently, the local solution
\(\zeta_*(C_\gamma)\) furnished by the implicit function theorem is the
unique positive critical point of \(h\) by (3.34), and hence the global maximum
point of \(h\).

Now, since \(\Phi_0(\zeta_*(C_\gamma)U_0)\le
\sup_{\zeta\ge0}\Phi_0(\zeta U_0)\), it follows from (3.31) and (3.36)
that
$$
h(\zeta_*(C_\gamma))
\le\sup_{\zeta\ge0}\Phi_0(\zeta U_0)
+\frac{C_\gamma}{\delta}\left[\zeta_*(C_\gamma)\right]^\delta
\le\sup_{\zeta\ge0}\Phi_0(\zeta U_0)
+\frac{C_\gamma}{\delta}(2\zeta_0)^\delta,
\qquad \forall C_\gamma<\varepsilon_1'.
\eqno(3.38)
$$
Set
$$
\varepsilon_1'':=\frac{\delta\eta_0}{2(2\zeta_0)^\delta}.
$$
Then
$$
\frac{C_\gamma}{\delta}(2\zeta_0)^\delta<\frac{\eta_0}{2},
\qquad \forall C_\gamma<\varepsilon_1''.
\eqno(3.39)
$$
Let \(\varepsilon_1:=\min\{\varepsilon_1',\varepsilon_1''\}\). For
\(C_\gamma<\varepsilon_1\), by (3.38), (3.39), (3.27) and (3.28),
$$
h(\zeta_*(C_\gamma))
<\sup_{\zeta\ge0}\Phi_0(\zeta U_0)+\frac{\eta_0}{2}
<\sup_{\zeta\ge0}\Phi_0(\zeta U_0)+\eta_0
=2^{-\frac{N-2s}{2s}}\frac{s}{N}S_K^{N/(2s)}.
$$
Therefore, it follows from (3.32) that
$$
\sup_{\zeta\ge0}\Phi(\zeta U_0)
\le \sup_{\zeta\ge0}h(\zeta)
=h(\zeta_*(C_\gamma))
<2^{-\frac{N-2s}{2s}}\frac{s}{N}S_K^{N/(2s)}.
$$

Hence \((S_1)\) holds with \(U_0\). Applying Theorem 1.1, the system
(1.2) admits a ground state solution
\(\widetilde U\in X\setminus\{0\}\) with
\(\Phi(\widetilde U)=\widetilde c>0\). \(\Box\)

\medskip
\noindent\textbf{Proof of Theorem 1.4.}
In the particular case \(K_i(x)=|x|^{-(N+2s)}\), the operator
\(\mathcal L_{K_i}\) coincides, up to a normalization constant, with the
fractional Laplacian \((-\Delta)^s\), and
\[
\widetilde E_i=\{v\in H^s(\mathbb R^N):\ v=0\ \text{a.e. in }
\mathbb R^N\setminus\Omega\},\qquad i=1,\dots,k,
\]
by Lemma 7-b in \cite{32}. Consequently,
\[
X=\big[\widetilde E_1\big]^k
=\{U=(u_1,\dots,u_k)\in [H^s(\mathbb R^N)]^k:\
u_i=0\ \text{a.e. in }\mathbb R^N\setminus\Omega,\ i=1,\dots,k\},
\]
and all the assumptions \((K_1)\)--\((K_3)\) are satisfied with
\(\rho=1\).

In this setting, the first eigenvalue and the best Sobolev constant of
the scalar problem become \(\sigma_{1,s}\) and \(S_s\), respectively,
the first eigenfunction becomes \(u_{0,s}\), and \((M_1)\) holds with
\(\sigma_1^i\) replaced by \(\sigma_{1,s}\), \(i=1,\dots,k\).

For cases (i) and (ii), applying Theorem 1.2 with
\(K_i(x)=|x|^{-(N+2s)}\), we obtain \((S_1)\). For case (iii), applying Theorem 1.3 with the same identification, and
taking \(\varepsilon_{1,s}\) to be the constant \(\varepsilon_1\)
furnished by Theorem 1.3 in this special case, we deduce from the
smallness condition
\[
\int_\Omega\gamma(x)|u_{0,s}|^\delta\,dx<\varepsilon_{1,s}
\]
that \((S_1)\) holds.

Thus, in all three cases (i)--(iii), condition
\((S_1)\) holds. Theorem 1.1 then yields a ground state solution
\(\widetilde U\in X\setminus\{0\}\) of system (1.4) satisfying
\(\Phi(\widetilde U)=\widetilde c>0\). \(\Box\)

\section*{Appendix}

\noindent\textbf{(Theorem 1.1, \cite{Jeanjean})}. {\it Let X be a Banach space equipped with the norm \(\|\cdot\|\) and let \(J\subset\mathbb{R}^+\) be an interval. We consider a family \(\{\Phi_\lambda\}_{\lambda\in J}\) of \(C^{1}\)-functionals on X of the form
$$
\Phi_\lambda(u) = A(u)-\lambda B(u),
$$
where \(B(u)\geq0,\ \forall u\in X\) and such that either \(A(u)\rightarrow+\infty\) or \(B(u)\rightarrow+\infty\) as \(\|u\|\rightarrow+\infty\). We assume that there are two points \((v_{1}, v_{2})\) in X such that setting
$$
\Gamma=\{\gamma\in C([0,1],X),\ \gamma(0)=v_1,\ \gamma(1)=v_2\}
$$
there hold, \(\forall\lambda\in J\),
$$
c_\lambda:=\inf_{\gamma\in\Gamma}\max_{t\in[0,1]}\Phi_{\lambda}(\gamma(t))>
\max\{\Phi_{\lambda}(v_1),\Phi_{\lambda}(v_2)\}.
$$
Then, for almost every \(\lambda\in J\), there is a sequence \(\{v_j\}\subset X\) such that \(\{v_j\}\) is bounded, \(\Phi_{\lambda}(v_j)\rightarrow c_\lambda\) and \(\Phi'_{\lambda}(v_j)\rightarrow0\) in the dual \(X^{-1}\) of \(X.\)}

\section*{Statements and Declarations}
This paper has no associated data.
The author declares that he has no
conflict of interest.

\end{document}